
\magnification 1200
\input mssymb.tex


\catcode`\@=11
\def\frak{\ifmmode\let\next\frak@\else
 \def\next{\errmessage{Use \string\frak\space only in math mode}}\fi\next}
\def\goth{\ifmmode\let\next\frak@\else
 \def\next{\errmessage{Use \string\goth\space only in math mode}}\fi\next}
\def\frak@#1{{\frak@@{#1}}}
\def\frak@@#1{\fam\euffam#1}
\def\Bbb{\ifmmode\let\next\Bbb@\else
 \def\next{\errmessage{Use \string\Bbb\space only in math mode}}\fi\next}
\def\Bbb@#1{{\Bbb@@{#1}}}
\def\Bbb@@#1{\fam\msyfam#1}
\catcode `\@=12

\font\smalltype=cmr10 at 10truept
 at 10truept
\nopagenumbers

\def\leftheadline{\ifnum\pageno > 1
		\tenrm\folio\hfil
		\smalltype Bertram, Daskalopoulos, and Wentworth
		\hfil\fi}

\def\rightheadline{\ifnum\pageno > 1
		\hfil\smalltype  Gromov Invariants
		\hfil\tenrm\folio\fi}

\headline={\ifodd\pageno\rightheadline\else\leftheadline\fi}

$\hbox{}$
\vskip .5in
\centerline{\bf GROMOV INVARIANTS FOR HOLOMORPHIC MAPS}
\smallskip
\centerline{\bf FROM RIEMANN SURFACES TO GRASSMANNIANS}
\vskip .75in
\centerline{Aaron Bertram
\footnote{$^1$}{\smalltype Supported in part by NSF Grant
DMS-9218215.}}
\smallskip
\centerline{Georgios Daskalopoulos}
\smallskip
\centerline{Richard Wentworth
\footnote{$^2$}{\smalltype Supported in part by NSF Mathematics
Postdoctoral Fellowship DMS-9007255.}}
\medskip
\centerline{\it Dedicated to Professor Raoul Bott on the occasion of
his 70th birthday}

\bigskip
\midinsert\narrower\narrower
\noindent {\bf Abstract.}
Two compactifications of the space of holomorphic
maps of  fixed degree from a compact Riemann surface to a
Grassmannian are studied.  It is shown that the Uhlenbeck
compactification has the structure of a projective variety and is
dominated by the algebraic compactification coming from the
Grothendieck Quot Scheme. The latter may be embedded into the moduli
space of solutions to a generalized version of the vortex
equations studied by Bradlow.  This gives an effective way of
computing certain intersection numbers (known as ``Gromov
invariants") on the space of holomorphic maps
into Grassmannians.  We carry out these
computations in the case where the Riemann surface has genus one.
\endinsert

\eject

\def\A{{\cal A}}
\def\D{{\frak D}}
\def\H{{\cal H}}
\def\F{{\cal F}}
\def\L{{\cal L}}
\def\M{{\cal M}}
\def\N{{\cal N}}
\def\U{{\overline{\M}_U}}

\def\Q{{\overline{\M}_Q}}
\def\QS{{\overline{\M}_Q^s}}
\def\Stau{{\cal V}_\tau}

\def\B{{\cal B}}
\def\O{{\cal O}}
\def\C{{\Bbb C}}
\def\P{{\Bbb P}}
\def\uk{{(\rho_\ast U)^{\oplus k}}}
\def\vk{{(\rho_\ast V)^{\oplus k}}}
\def\uktilde{{(\tilde\rho_\ast \tilde U)^{\oplus k}}}
\def\vktilde{{(\tilde\rho_\ast \tilde V)^{\oplus k}}}
\def\pspace{{\overline{\M}_P}}
\def\pspaceo{{\overline{\M}_P^o}}

\def\data{(d,r,k)}
\def\tdata{(\tilde d,r,k)}
\def\xtuple{(X_1,\ldots, X_r)}
\def\ztuple{(Z_1,\ldots, Z_r)}
\def\rk{\hbox{rk}\,}
\def\G{{\frak G}}
\def\Go{{\frak G}_0}
\def\Gc{{\frak G}^\C}
\def\Gcbar{\overline{\frak G}^\C}

\def\Gcp{{\frak G}^\C_p}
\def\dbar{{\bar\partial}}
\def\dbare{{\dbar_E}}
\def\kpair{{(\dbare,\vec\phi)}}
\def\curv{F_{\dbare,H}}
\def\cplx{C^\dbare_\phi}
\def\endo{\hbox{End}\,}
\def\lie{\hbox{Lie}\,}
\def\tr{\hbox{tr}\,}

\def\I{I(d,k; n)}

\def\lra{\longrightarrow}
\def\mapright#1{\,\smash{\mathop{\longrightarrow}
	\limits^{#1}}\,}
\def\mapdown#1{\downarrow\rlap
	{$\vcenter{\hbox{$\scriptstyle#1$}}$}}

\baselineskip 16pt \normalbaselineskip 16pt

\beginsection 1. Introduction

In [G1] Gromov introduced the moduli space of pseudo-holomorphic
curves in order to obtain new global invariants of symplectic
manifolds.  The theory has turned out to have far-ranging
applications, perhaps the most spectacular being to the Arnold
conjecture and the subsequent development of Floer homology (see
[F]).  The Donaldson-type version of Gromov's invariants [G2] have
not been much studied by mathematicians (see, however, [R]),
although in the physics literature they arise
naturally in the topological quantum field theories that
have attracted so much attention recently (see [Wi]).  In this paper
we shall provide a framework for the calculation of these invariants
for the case of holomorphic maps from a fixed compact Riemann
surface $C$ of genus $g$ to the Grassmannian $G(r,k)$ of complex
$r$-planes in $\C^k$.

To introduce some notation, let $\M$ (or $\M\data$) denote the
space of holomorphic maps $f:C\to G(r,k)$ of degree $d$.  For
simplicity we shall usually omit the data $d,r$ and $k$, and
since $C$ will always be fixed we suppress it from the notation as
well.
Furthermore, we will assume throughout the paper, unless otherwise
indicated, that $g\geq 1$ and $d > 2r(g-1)$.  The evaluation map
$$\mu:C\times\M\longrightarrow G(r,k): (p,f)\mapsto
f(p)\leqno(1.1)$$
defines cohomology classes on $\M$ by pulling back classes from
$G(r,k)$ and slanting with the homology of $C$ (we will always slant
with a point).  Thus classes $X_1,\ldots, X_l\in
H^\ast\left(G(r,k),\C\right)$ of dimensions $n_1,\ldots, n_l$
and  a point $p\in C$ define
classes $\mu^\ast X_1/p,\ldots,$ \break
$\mu^\ast X_l/p
\in H^\ast(\M,\C)$, and the invariants we are interested in are the
intersection numbers
$$\langle X_1^{s_1}\cdots X_l^{s_l}\rangle =:
(\mu^\ast X_1/p)^{s_1}\cup\cdots\cup(\mu^\ast X_l/p)^{s_l}[\M]\;,
\leqno(1.2)$$
where $\sum s_i n_i$ is the dimension of $\M$.  As in Donaldson
theory there is the usual difficulty in making sense of (1.2); one
of the first problems is to find a suitable compactification of
$\M$.  This leads us to three different approaches to the problem
which may roughly be characterized as the analytic, algebraic, and
gauge theoretic points of view, as we shall now explain.

The analytic approach compactifies $\M$ through the mechanism of
``bubbling" presented in the fundamental work of Sacks and Uhlenbeck
[S-U].  By adding to $\M$ the data of a divisor on $C$ and a
holomorphic map of lower degree, one obtains a compact topological
space $\U$ which contains $\M$ as an open set (in [G1] the
compactness of $\U$ is proved by different methods; see [Wf]
for the details of applying the results of [S-U] to this situation).
We shall refer to $\U$ as the {\it Uhlenbeck compactification}
of $\M$.  In
\S 4 we shall prove

\proclaim Theorem 1.3.  $\U\data$ has the structure of a projective
variety.  For the case of maps to projective space (i.e. $r=1$),
$\U(d,1,k)$ is smooth and is in fact a projective bundle over $J_d$, the
Jacobian variety of degree $d$ line bundles on $C$.

Taking a cue from the theory of four manifolds, we may expect there
to be an algebro-geometric compactification of $\M$ which perhaps
contains more information than $\U$.  This is indeed the case, and
the appropriate object is a particular case of the {\it Grothendieck
Quot Scheme}, which we denote by $\Q\data$.
This is a projective variety
parameterizing  quotients
$$\O_C^k\longrightarrow {\cal F}\longrightarrow 0\;,$$
where $\O_C$ is the structure sheaf of $C$ and $\cal F$ is a
coherent  sheaf on $C$ with a given Hilbert polynomial determined by
$r$ and $d$. It will be seen that $\Q$  naturally contains $\M$ as
an open subvariety.  The relationship with the Uhlenbeck
compactification is given by the following

\proclaim Theorem 1.4.  There is an algebraic surjection
$$\Q\data\mapright{u} \U\data$$
which is an isomorphism on $\M$.  For the case of maps to projective
space, $u$ itself is an isomorphism.
Moreover, for $d$ sufficiently large (depending on $g$, $r$, and $k$)
$\M\data$ is dense, and
the scheme structures on both $\Q\data$ and $\U\data$ are
irreducible and generically reduced.

The map $u$ gives us a way to lift the calculation of intersections
on $\U$ to $\Q$.  However, in general the variety $\Q$  is still
singular, so it would be nice to have an explicit embedding of $\Q$
in a smooth variety where the intersecting classes extended.  This
leads us to the third, gauge theoretic construction:

We extend Bradlow's notion of {\it stable pairs} to the case of
a rank $r$ holomorphic bundle $E\to C$ with $k$ holomorphic
sections $\phi_1,\ldots, \phi_k$ (cf. [B]).  The definition of
stability will depend on the choice of a real parameter $\tau$, and
we shall see that the $(k+1)$-tuple $(E; \phi_1,\ldots,\phi_k)$ is
$\tau$-semistable if and only if it admits a solution to a certain
non-linear PDE which we call the $k$-$\tau$-{\it vortex equation} (see
Theorem 3.5).  Then following [B-D], we construct the moduli space
$\B_\tau$ of solutions.  We prove

\proclaim Theorem 1.5.  For generic $\tau$ in a certain
admissible range (see Assumption 3.12),\break
$\B_\tau\data$ is a smooth, projective variety.  For $r=1$ and
$\tau > d$,
$$\B_\tau(d,1,k)\simeq\Q(d,1,k)\simeq\U(d,1,k)$$
as projective varieties.

In the process of constructing
$\B_\tau$ we also verify certain universal
properties and the existence of a universal
rank $r$ bundle $U_\tau\to C\times\B_\tau$.  This will allow
us to prove

\proclaim Theorem 1.6.  For a given $d$, there exists a choice of
$\tilde d > d$ and a choice of $\tau$ within the range of Theorem 1.5
such that we have an algebraic embedding
$$\Q\data\mapright{b} \B_\tau\tdata\;.$$
Moreover, the image $b(\Q)$ is cut out by equations determined by the
top Chern class of $U_\tau\tdata$.

\noindent
The explicit equations are the strength of Theorem 1.6;
they are, of course, necessary if we want to
compute intersection numbers. Other equations
arise if
we simply map $\U\data\to\U(d,1,{k\choose r})$ via the Pl\"ucker
embedding, as we shall see in \S 4; however,
as with the Pl\" ucker  embedding itself, there are too many
equations, and from the point of view of computing intersections
this is therefore not useful.  By contrast, the equations in
Theorem 1.6  cut out the Quot scheme as a complete intersection,
least for $d$ sufficiently large.

The final step in our algorithm is to extend the classes to
$\B_\tau\tdata$ and reduce the computation  of intersection numbers
to $\N(\tilde d, r)$, the Seshadri moduli space of rank $r$ semistable
bundles on $C$ of degree $\tilde d$.  The structure of
$H^\ast(\N(\tilde d, r),\C)$ is understood, at least in principle,
from the results of Kirwan [K], or for $r=2$ and $\tilde d$ odd,
from [T1] and [Z].

Actually, at present the reduction to $\N$ only works for
$r=2$.  The idea is to study the behavior of $\B_\tau$ as
the parameter $\tau$ passes through a non-generic value, in a manner
similar to the description given by Thaddeus
in [T2] (see also [B-D-W]).  The
relationship is as follows:  If $\tau$ is a non-generic value and
$\varepsilon$ is small, then $\B_{\tau-\varepsilon}$ is gotten
from $\B_{\tau+\varepsilon}$
by blowing up along a smooth subvariety and then blowing down in a
different direction; a process known as a ``flip".  This structure
arises naturally from realizing the parameter $\tau$ as the Morse
function associated to a circle action on some bigger symplectic
manifold and then applying the  results of Guillemin and Sternberg
[G-S].  In this way one eventually arrives at a value of $\tau$ for
which $\B_\tau$ is a projective bundle over $\N$
(for even $\tilde d$, the fiber is only generically a projective
space).  One can keep track of how the intersection numbers change
as $\B_\tau$ is flipped, and this in principle reduces
the calculation to intersections on $\N$.

The computations are, of course, rather unwieldy, but to demonstrate
that this procedure can actually be carried through we compute in
\S 5.3 the case of maps from an
elliptic curve to $G(2,k)$. In this case, $\dim\Q(d,2,k)=kd$ for $d$
large.  The intersection numbers $\langle X_1^{kd-2n}X_2^n\rangle$,
$n=0,1,\ldots,[kd/2]$ are rigorously defined in \S 5.1.  We then
prove

\proclaim Theorem 1.7.  For holomorphic maps of degree $d$
sufficiently large from an elliptic curve to $G(2,k)$, the
intersection numbers are given by
$$\langle X_1^{kd-2n}X_2^n\rangle
=(-1)^{d+1}k 2^{kd-2n-1}-
(-1)^{d+1}{k^2\over 2}
\sum_{p\in{\Bbb Z}\atop n/k\leq p\leq d-n/k}{kd-2n\choose
kp-n} \;.$$

In \S 5.2, we briefly discuss a
remarkable conjecture concerning all the intersection  numbers on
$\U$ which is due to Vafa and Intriligator (see [V], [I]).  Using
arguments from the physics of topological sigma models they derive a
formula for the numbers which is based entirely on a residue
calculation involving the homogeneous polynomial which characterizes
the cohomology
ring $H^\ast\left( G(r,k),\C\right)$.  Somewhat surprisingly,
this formula, which is simple to state yet is highly non-trivial
(see Conjecture 5.10), agrees with
our result Theorem 1.7 above.

As a final note, the case $r=1$, i.e. maps into projective space,
has a simple description.  The intersection numbers can all be
computed, and we will do so in \S 2.  It is found that these
also agree with the physics formula (see Theorems 2.9 and 5.11).

\beginsection 2.  Maps to projective space

The purpose of this section is to compute the Gromov invariants
$\langle X^m\rangle$ of (1.2) for the case of holomorphic maps into
projective spaces.  It turns out that in this case there is an
obvious nonsingular compactification of
the space of holomorphic maps.  This
can be described in terms of the push-forward of the universal
line bundle on the Jacobian variety, and the computation of the
intersection numbers reduces to the Poincar\'e formula for the theta
divisor.

Let $C$ be a compact Riemann surface of genus $g$
with a base point $p\in C$. {\it In this section
only} we include the case $g=0$.  Let
$\M=\M(d,1,k)$ denote the space of holomorphic maps from $C$ to
$\P^{k-1}$ of non-negative degree $d > 2g-2$.
It is well-known that $\M$ is a complex manifold of dimension
$m=kd-(k-1)(g-1)$.  Let $H$ denote a fixed hyperplane in
$\P^{k-1}$.  Then we define
$$X=\left\{ f\in\M : f(p)\in H \right\}\;.\leqno(2.1)$$
Clearly, $X$ is a divisor in $\M$.  Our goal is to compute the top
intersection $\langle X^m\rangle$, where $m=\dim\M$ as above.

Since $\M$ is not compact, in order to make sense of $\langle
X^m\rangle$ we shall define a smooth compactification
$\pspace$ of $\M$ and extend the class $X$ to $\pspace$.
Then we can define:
$$\langle X^m\rangle = X^m[\pspace]\;,\leqno(2.2)$$
where we denote both the extended class and its Poincar\'e dual by $X$,
and by $[\pspace]$ we mean the fundamental class of $\pspace$.

Let $J_d$ denote the
Jacobian variety of degree $d$ line bundles on $C$, and let $U\to
C\times  J_d$ denote the universal or Poincar\'e line bundle
on $C\times J_d$.
By this we mean a line bundle whose
 restriction  to  $C\times \{L\}$ is a line bundle on $C$
isomorphic to $L$ (cf. [A-C-G-H]).
This bundle is not uniquely determined, since we are free to tensor
a given choice with any line bundle on $J_d$.
It will be convenient to normalize $U$ so that its restriction
to the point $p\in C$:
$$U_p= U\bigr|_{\{p\}\times J_d}\leqno(2.3)$$
is holomorphically trivial.
Let $\rho : C\times J_d\to J_d$ be the projection
map.  By our requirement that $d> 2g-2$,
$\rho_\ast U$ is a vector bundle on $J_d$  of rank $d-g+1$.  Now
define
$$\pspace=\pspace(d,k)
=\P\left((\rho_\ast U)^{\oplus k}\right)\;,\leqno(2.4)$$
where the superscript $\oplus k$ means the
fiberwise direct sum of $k$ copies
of the vector bundle $\rho_\ast U$.
Thus $\pspace\mapright{\pi} J_d$ is a projective bundle with fiber over
$L$ isomorphic to the projective space $\P\left(
H^0(C,L)^{\oplus k}\right)$.

Alternatively, $\pspace$ may be thought of as
gauge equivalence classes of what we shall henceforth refer to
as {\it k-pairs}, by which we mean  $(k+1)$-tuples
$(L; \phi_1,\ldots,\phi_k)$
where the $\phi_i$'s are holomorphic sections of $L$ and
$\vec\phi =: (\phi_1,\ldots,\phi_k)\not\equiv(0,\ldots,0)$.  The latter
condition implies that the sections generate the fiber of $L$ at a
generic point in $C$.

In order to show that $\pspace$ forms a compactification of $\M$,
let $\pspaceo\subset \pspace$ denote the open subvariety
consisting of  those $k$-pairs for which the set of $\phi_i$'s
generates the fiber at {\it every} point. Then we define a map
$$F: \pspaceo\longrightarrow \M\leqno(2.5)$$
as follows:  Given a point $[L,\vec\phi]\in\pspaceo$,
let $(L,\vec\phi)$ be a representative, and let $f:
C\to\P^{k-1}$ be the map
$$f(p)=[\phi_1(p),\ldots,\phi_k(p)]\;.\leqno(2.6)$$
Since $\vec\phi(p)\neq 0$ for every $p$, this is a well-defined
holomorphic map which is easily seen to have degree $d$.
Since a different choice of
representative of $[L,\vec\phi]$ has the effect
of rescaling each $\phi_i(p)$ by the same  constant  $\lambda(p)\in
\C^\ast$, the map $f$ is independent of this choice, and so $F$ is
well-defined and clearly holomorphic.
Conversely, given $f\in\M$, let $S$ denote the tautological line
bundle on $\P^{k-1}$.  Then $L=f^\ast S^\ast$ is a line bundle of
degree $d$ on $C$, and the coordinates of $\C^k$ pull-back to
$k$ holomorphic sections $\phi_1,\ldots, \phi_k$ of $L$
generating the fiber at every point.  Clearly, the
map $f\mapsto [L,\vec\phi]$ is a holomorphic inverse of $F$.
Therefore, $F$ is a biholomorphism and we have

\proclaim Proposition 2.7.  The projective bundle
$\pspace\mapright{\pi} J_d$ is a compactification of $\M$.

\noindent
We will see in \S 4 that $\pspace$ coincides with the
Uhlenbeck compactification of $\M$.

Now by the top intersection of $X$ on $\M$ we shall mean
the top intersection of its Zariski closure in $\pspace$.
Let us denote this extension and its Poincar\'e dual also by $X$.
In order to compute the intersection number (2.2) we proceed as
follows:
Let $\O_\pspace(1)$
denote the anti-tautological line bundle on $\pspace$ and
$U_p\to J_d$ the restriction of $U$ to $\{p\}\times J_d$.
Then we have

\proclaim Lemma 2.8.  $c_1(\O_\pspace(1))= X$.

\noindent {\it Proof.}  Given a line bundle $L\to C$,
consider the map
$$\psi_p: H^0(C,L)\times\cdots\times H^0(C,L)\longrightarrow
L_p: (\phi_1,\ldots,\phi_k)\longmapsto \phi_1(p)\;.$$
Then since $U_p$ is trivial,
$\psi_p$ is a well-defined
linear form on the homogeneous coordinates of the
fiber of $\pspace\to J_d$.  It follows that
$$s\left([L,\vec\phi]\right)= (\ker\psi_p)^\ast$$
defines a holomorphic section of $\O_\pspace(1)$,
and the zero locus of $s$ is
$$Z(s)=\left\{ [L,\vec\phi] : \phi_1(p)=0 \right\}\;.$$
The isomorphism $F$ in (2.5) identifies $Z(s)\cap \pspaceo$ with
the subspace $X$ in (2.1) (for $H$ the hyperplane defined by
$[0,z_2,\ldots, z_k]\in \P^{k-1}$), and
this completes the proof of the lemma.

Now we are prepared to prove the main result of this section.

\proclaim Theorem 2.9.  For $X$ defined as in (2.1), non-negative $d >
2g-2$, and $m=\dim\pspace=kd-(k-1)(g-1)$  we have
$\langle X^m\rangle =k^g$.

The computation is based on standard results on the cohomology ring
of projectivized bundles.
The most important tool is the notion of a Segre class; since this
is perhaps not so well-known to analysts, we briefly review the
essentials.

Let $V$ be a rank $r$ holomorphic vector bundle on a compact,
complex manifold $M$ of dimension $n$, and let $\P(V)$ denote the
projectivization of $V$.  Then $\P(V)\mapright{\pi} M$ is a
projective bundle with fiber $\P^{r-1}$.
We then have an exact sequence of bundles
$$0\lra\O_\P(-1)\lra\pi^\ast V\lra Q \lra 0\;,$$
where $\O_\P(-1)$ denotes the tautological line bundle on $\P(V)$ and
$Q$ the quotient rank $r-1$ bundle.
Let $X=c_1(\O_\P(1))$. Then
$$(1-X)c(Q)= \pi^\ast c(V)\;,$$
where $c$ denotes the total Chern polynomial, or equivalently
$$c(Q)=\pi^\ast c(V) (1+X+X^2+\cdots)\;.$$
By applying the push-forward homomorphism $\pi_\ast$, or integration
along the fibers, and using the fact that $\pi_\ast c_i(Q)=0$ for
$i< r-1$ and that $c_{r-1}(Q)$ restricts to the fundamental
class of the fiber, we obtain
$$1= c(V)\pi_\ast(1+X+X^2+\cdots)\;.\leqno(2.10)$$
We now define the {\it total Segre class} of $V$
by the formal expansion
$$s(V) = {1\over c(V)}\;,\leqno(2.11)$$
and the {\it Segre classes} $s_i(V)$  of $V$ are defined to be the
$i$-th homogeneous part of $s(V)$. It follows from (2.10) and (2.11)
that for every $l\geq r-1$,
$$\pi_\ast X^l = s_{l-r+1}(V)\;.\leqno(2.12)$$

\noindent {\it Proof of Theorem 2.9.}
The discussion above applies to our situation if we
let $V=(\rho_\ast U)^{\oplus k}$.
This is a vector bundle of rank $k(d+1-g)$ on $J_d$, and
$$\pspace = \P(V)\mapright{\pi} J_d\;.$$
Moreover, by Lemma 2.8  and (2.12),
$$\langle X^m\rangle= X^m[\pspace]=\pi_\ast
X^{k(d+1-g)-1+g}[J_d]=s_g(V)[J_d]\;.$$
It therefore suffices to compute $s_g(V)$.
The Chern character of $V$ is computed by the
Grothendieck-Riemann-Roch formula (see [A-C-G-H], p. 336)
$$ch\left((\rho_\ast U)^{\oplus k}\right)=k\cdot
ch(\rho_\ast U)= k(d-g+1)-k\theta\;,$$
where $\theta$ denotes the dual of the theta divisor in $J_d$.
Moreover, as in [A-C-G-H], p. 336,
the expression above implies a particularly nice form
for the Chern polynomial, $c(V)=
e^{-k\theta}$, and hence $s_g(V)=k^g\theta^g/g!$.
Applying the Poincar\'e formula (see [A-C-G-H], p. 25)
completes the proof.

It is somewhat curious that this is precisely the dimension of the
space of level $k$ theta functions for genus $g$.  We shall see in
\S 5.2 that Theorem 2.9 confirms the physics conjecture for maps
to projective space.

\beginsection 3.  Moduli of stable k-pairs

\bigskip
\centerline{\it \S 3.1 Definition of stability}

\noindent In this section we generalize the notion of stable pairs
to stable $k$-pairs.  We give the precise definition of stability
for $k$-pairs and describe the associated Hermitian-Einstein
equations.  Since most of this section is a direct generalization of
the corresponding results for stable pairs, we shall give only a
brief exposition and refer to [B], [B-D] and [Ti] for further
details.

Let $C$ be a compact Riemann surface of genus $g\geq 1$ and $E$ a
complex vector bundle on $C$ of rank $r$ and degree $d > 2r(g-1)$.
Unless otherwise stated, we assume that we have a fixed K\"ahler
metric on $C$ of area $4\pi$ and a fixed hermitian metric on $E$.

Let $\D$ denote the space of $\dbar$-operators on $E$, and let
$\Omega^0(E)$ denote the space of smooth sections of $E$.  We
topologize both $\D$ and $\Omega^0(E)$ by introducing the
appropriate Sobolev norms as in [B-D].
The space of $k$-{\it pairs} is defined to be
$$\A=\A\data=\D\times\Omega^0(E)\times\cdots\times\Omega^0(E)\;,\leqno
(3.1)$$
where we take $k$ copies of $\Omega^0(E)$.  For example, $\A(d,r,1)$
is the space of {\it pairs} considered in [B-D].  The space of {\it
holomorphic k-pairs} is defined to be
$$\eqalign{\H=\H\data &=\biggl\{ (\dbare,\phi_1,\ldots,\phi_k)\in\A\data
: \dbare\phi_i=0, i=1,\ldots,k\;,\cr
&\qquad\qquad\qquad \hbox{ and }(\phi_1,\ldots,\phi_k)\not\equiv
(0,\ldots,0)\biggr\}\;.\cr}\leqno(3.2)$$
We shall often denote the $k$-pair $(\dbare,\phi_1,\ldots,\phi_k)$
by  $(\dbare,\vec\phi)$.
To introduce the notion of stability, define as in [B], [Ti], the
numbers
$$\eqalign{
\mu_M(E)&=\hbox{max}\biggl\{ \mu(F) : F\subset E\hbox{ a holomorphic
subbundle with } \rk(F) >0\biggr\}\cr
\mu_m(\vec\phi)&=\hbox{min}\biggl\{ \mu(E/E_\phi) : E_\phi
\subset E\hbox{ a proper holomorphic
subbundle with}\cr
&\qquad\qquad\qquad\qquad
\phi_i\in H^0(E_\phi), i=1,\ldots,k\biggr\}\;.\cr}$$
Here $\mu$ denotes the usual Schatz slope $\mu=\deg/\rk$.  Note that
in the definition of $\mu_m(\vec\phi)$ the set of such
$E_\phi$'s may be empty, in which case we set
$\mu_m(\vec\phi)=+\infty$.  Of course,
if the rank is two or greater this cannot happen with
pairs, i.e. $k=1$.

\proclaim Definition 3.3.   A holomorphic $k$-pair
$(\dbare,\vec\phi)\in\H$ is called $\tau$-{\it stable} for $\tau\in{\Bbb
R}$ if
$$\mu_M(E) < \tau < \mu_m(\vec\phi)\;.$$
A $k$-pair is called {\it stable} if it is $\tau$-stable for some
$\tau$.

Note that if $k=1$ this definition agrees with the one in [B].  We
shall denote by $\Stau=\Stau\data$ the subspace of $\H\data$
consisting of $\tau$-stable $k$-pairs.

It is by now standard philosophy that any reasonable stability
condition corresponds to the existence of special bundle metrics.
These metrics satisfy the analogue of the Hermitian-Einstein
equations, which we now describe.
Given an hermitian metric $H$ on $E$ and $\phi$ a smooth section of
$E$, we denote by $\phi^\ast$ the section of $E^\ast$ obtained by
taking  the hermitian adjoint of $\phi$.  Also, given a
$\dbar$-operator $\dbare$ on $E$ we denote by $\curv$
the curvature of the unique hermitian connection compatible with
$\dbare$ and $H$.  Finally, we define the $k$-$\tau$-{\it
vortex equation} by
$$\sqrt{-1}\ast\curv+{1\over
2}\sum_{i=1}^k\phi_i\otimes\phi_i^\ast={\tau\over 2}{\bf I}\;.
\leqno(3.4)$$
We now have the following

\proclaim Theorem 3.5. Let $\kpair\in\H$ be a holomorphic $k$-pair.
Suppose that for a given value of the parameter $\tau$ there is an
hermitian metric $H$ on $E$ such that the $k$-$\tau$-vortex
(3.4) is satisfied.  Then $E$ splits holomorphically $E=E_\phi\oplus
E_s$, where
\itemitem{(i)} $E_s$, if nonempty, is a direct sum of stable bundles
each of slope $\tau$.
\itemitem{(ii)} $E_\phi$ contains the sections $\phi_i$ for all
$i=1,\ldots, k$, and with the induced holomorphic structure from $E$
the $k$-pair $(E_\phi,\vec\phi)$ is $\tau$-stable.

\noindent{\sl Conversely,
suppose that $\kpair$ is $\tau$-stable.  Then the
$k$-$\tau$-vortex equation has a unique solution.}

\noindent {\it Proof.}  See [B] and [Ti].

\bigskip
\centerline{\it \S 3.2 Moduli of  k-pairs and universal bundles}

\noindent In this section we construct for generic values of $\tau$
within a certain range a smooth moduli space $\B_\tau=\B_\tau\data$
of $\tau$-stable $k$-pairs on the vector bundle $E$.  Furthermore,
we construct a universal rank $r$ bundle $U_\tau \to C\times\B_\tau$
with $k$ universal sections.

We first recall that the natural complex structure
on $\D$ and $\Omega^0(E)$ induces a complex structure on the space
of $k$-pairs $\A$.  This in turn defines a complex structure on the
space of holomorphic $k$-pairs $\H\subset\A$.  In other words, $\H$
has the structure of an infinite dimensional analytic variety.  It
was observed in [B-D-W], Corollary 2.7,
that if $k=1$, $\H$ is a complex manifold.
However, this fails to be true for general $k$, and for the purpose
of this paper we will restrict to an open smooth part $\H^\ast$ of
$\H$ defined as follows:

\proclaim Definition 3.6.  If $r=1$, let $\H^\ast=\H$.  If $r\geq
2$, let $\H^\ast$ consist of those $k$-pairs $\kpair\in\H$
satisfying
$$\mu_M(E) < {d-(2g-2)\over r-1}\;.$$

\proclaim Proposition 3.7.  $\H^\ast$ is a smooth, complex
submanifold of $\A$.

\noindent {\it Proof.}  Let
$G : \A\to \Omega^1(E)\times\cdots\times\Omega^1(E)$
be the map $G\kpair=(\dbare\phi_1,\ldots,\dbare\phi_k)$. Given
$\kpair\in\H$, the derivative of $G$ at $\kpair$ is given by (cf.
[B-D], (2.11))
$$\delta G_\kpair(\alpha,\eta_1,\ldots,\eta_k)=
(\dbare\eta_1+\alpha\phi_1,\ldots, \dbare\eta_k+\alpha\phi_k)\;.$$
If the cokernel of $\dbare$ vanishes for all $\kpair\in\H^\ast$,
then the restriction of $G$ to $\H^\ast$ has everywhere
derivative of maximal rank, and Proposition 3.7 then follows from
the implicit function theorem.  It therefore suffices to prove the
following

\proclaim Lemma 3.8.  Let $E$ be a holomorphic rank $r$
vector bundle on $C$ of degree $d$, $d> 2r(g-1)$.
Suppose that  either $r=1$ or $r>1$  and
$$\mu_M(E) < {d-(2g-2)\over r-1}\;.$$
Then $H^1(E)=0$.

\noindent {\it Proof.}  The case $r=1$ is  a vanishing theorem. Suppose
therefore that $r>1$.  Consider the Harder-Narasimhan filtration of
the holomorphic bundle $E$
$$0=E_0\subset E_1\subset\cdots\subset E_l=E\;,\leqno(3.9)$$
where $D_j=E_j/E_{j-1}$ is semistable and $\mu_j=\mu(D_j)$ satisfies
$\mu_1 >\cdots >\mu_l$.  In particular, we have exact sequences
$$0\lra E_{j-1}\lra E_j\lra D_j\lra 0\;,$$
where $j=1,\ldots, l$.
By induction, it suffices to show $H^1(D_j)=0$.  Since $D_j$ is
semistable and $\mu_1 >\cdots >\mu_l$, it is enough to show that
$\mu_l > 2g-2$ (cf. [New], p. 134).
This can be proved as follows: First
write
$$d=\deg E=\deg E_{l-1} +\deg D_l\;.\leqno(3.10)$$
By assumption,
$$\mu(E_{l-1})\leq\mu_M(E) <{d-(2g-2)\over r-1}\;.\leqno(3.11)$$
 From (3.10) and (3.11) we obtain
$$\eqalign{\deg D_l &> d-(r-\rk D_l)\left({d-(2g-2)\over r-1}\right)\cr
&={r(2g-2)-d\over r-1}+\rk D_l\left({d-(2g-2)\over
r-1}\right)\;.\cr}$$
Since $r(2g-2)-d < 0$, we obtain
$$\deg D_l > \rk D_l\left\{ {r(2g-2)-d\over r-1}+{d-(2g-2)\over
r-1}\right\}=\rk D_l (2g-2)\;.$$
Thus $\mu_l > 2g-2$, which completes the proof of Lemma 3.8 and
Proposition 3.7.

For our construction of the moduli space, we make the following

\proclaim Assumption 3.12. The {\it admissible range} of $\tau$ is
defined as
$${d\over r}< \tau <{d-(2g-2)\over r-1}\;.$$

\proclaim Definition 3.13.  A value of $\tau$
is called {\it generic} if $\tau$ is {\it not} rational of the form $
p/q$ where $0< q < r$.

The construction of a moduli space of $\tau$-stable $k$-pairs will
only apply for $\tau$ in the admissible range, however, we shall
still be interested in other values of $\tau$.  Specifically, we
note the following

\proclaim Proposition 3.14.  Given $\tau > d$, a holomorphic
$k$-pair $\kpair$ is $\tau$-stable if and only if the sections
$\{\phi_1,\ldots,\phi_k\}$ generically generate the fiber of $E$ on
$C$.  In particular, if $k< r$ then the range of values $\tau$ for
which there exist $\tau$-stable $k$-pairs is bounded.

\noindent {\it Proof.}  If $\{\phi_1,\ldots,\phi_k\}$ do not
generically generate the fiber of $E$, then they fail to generate at
every point.  Hence they span a proper subbundle $E_\phi\neq E$.
But then
$$\mu_m(\vec\phi)\leq \mu(E/E_\phi)\leq \deg (E/E_\phi)\leq d <
\tau\;,$$
and so $\kpair$ cannot be $\tau$-stable.  On the other hand, if
$\{\phi_1,\ldots,\phi_k\}$ generate the fiber of $E$ generically,
then they cannot span a proper subbundle, and so
$\mu_m(\vec\phi)=\infty$.  Thus to prove $\tau$-stability it
suffices to show that $\mu_M(E)\leq d$.  Consider once again the
filtration (3.9).  For any $E'\subset E$ we have
$\mu(E')\leq\mu(E_1)=\mu_1$ (see [New], p. 162).  Therefore, we need
only show $\deg E_1\leq d$.  We have an exact sequence
$$0\lra E_{l-1}\lra E\lra E/E_{l-1}\lra 0\;.$$
Now $E/E_{l-1}$ is semistable and has non-trivial sections, since
the fiber of $E$ is supposed to be generated generically by
$\{\phi_1,\ldots,\phi_k\}$.  Thus, $\deg E/E_{l-1} \geq 0$.
Inductively, suppose that $\deg E/E_i\geq 0$.  Then we have
$$0\lra E_i/E_{i-1}\lra E/E_{i-1}\lra E/E_i\lra 0\;.$$
We have shown that $\mu_l\geq 0$, so in particular $\deg E_i/E_{i-1}
> 0$ for $i< l$.  Hence,
$$\deg E/E_{i-1}=\deg E/E_i +\deg E_i/E_{i-1}\;,$$
which implies $\deg E/E_{i-1}\geq 0$.  Finally, consider
$$0\lra E_1\lra E\lra E/E_1\lra 0\;.$$
Then $\deg E_1=\deg E -\deg E/E_1\leq d$.
This completes the proof.

After this digression,
we are now ready to proceed with the construction of the moduli
space. Recall from \S 3.1 the subspace $\Stau\subset\H$ of
$\tau$-stable $k$-pairs.  For $\tau$ admissible in the sense of
Assumption 3.12, $\Stau$ is an open
submanifold  of $\H^\ast$ and is therefore a smooth manifold.
Moreover, the actions of the complex gauge group on $\D$ and
$\Omega^0(E)$ give an action on $\Stau$.  We define
$\B_\tau=\B_\tau\data$ to be the quotient of $\Stau$ by $\Gc$.
Our next task is to put a complex manifold structure on $\B_\tau$.
Following [B-D], for $\kpair\in\H$ we define the
complex
$$\cplx : \Omega^0(\endo E)\mapright{d_1}\Omega^{0,1}(\endo
E)\oplus\Omega^0(E)\oplus\cdots\oplus\Omega^0(E)\mapright{d_2}
\Omega^{0,1}(E)\oplus\cdots\oplus\Omega^{0,1}(E)\;,$$
where
$$\eqalign{
d_1(u)&=(-\dbare u, u\phi_1,\ldots, u\phi_k)\cr
d_2(\alpha,\eta_1,\ldots,\eta_k)&=(\dbare\eta_1+\alpha\phi_1,\ldots,
\dbare\eta_k+\alpha\phi_k)
}\leqno(3.15)$$
The properties of $\cplx$ that we need may be summarized as follows:

\proclaim Proposition 3.16.  Let $\kpair\in \H$.  Then
\itemitem{(i)} $\cplx$ is an elliptic complex.
\itemitem{(ii)} If $\kpair\in\H^\ast$, then $H^2(\cplx)=0$.
\itemitem{(iii)} If $\kpair$ is $\tau$-stable, then $H^0(\cplx)=0$.
\itemitem{(iv)} For $\kpair\in\H^\ast$,
$\chi(\cplx)= kd  - r(k-r)(g-1)\;.$

\noindent {\it Proof.}
Part (i) follows as in [B-D], Proposition 2.1, and
part (ii) follows immediately from Lemma 3.8.  Part
(iii) is proved as in [B-D], Proposition 2.3.  Finally, for part
(iv), observe that as in [B-D], Proposition 3.6,
$$\chi(\cplx)=k\chi(E)-\chi(\endo E)= kd-r(k-r)(g-1)\;.$$

A complex slice theorem as in [B-D], \S 3 now proves

\proclaim Proposition 3.17.  For $\tau$ satisfying Assumption 3.11,
$\B_\tau $ is a complex manifold of dimension $kd-r(k-r)(g-1)$.
Moreover, its tangent space may be identified
$$T_{[\dbare,\vec\phi]}\B_\tau=H^1(\cplx)\;.$$

We next define a K\"ahler structure on $\B_\tau$.  Recall that $\D$
and $\Omega^0(E)$ have natural K\"ahler forms $\Omega_\D$ and
$\Omega_{\Omega^0(E)}$, compatible with the $L^2$-inner products
(cf. [B-D], \S 4).  More precisely, let
$$\eqalign{
\Omega_\D(\alpha,\beta)&=\sqrt{-1}\left(\langle\alpha,\beta\rangle_\D-
\langle\beta,\alpha\rangle_\D\right)\cr
\Omega_{\Omega^0(E)}(\eta,\nu)&={\sqrt{-1}\over
2}\left(\langle\eta,\nu\rangle_{\Omega^0(E)}-\langle\nu,\eta\rangle_{
\Omega^0(E)}\right)\;.\cr}$$
These combine to define a K\"ahler form
$$\Omega=\Omega_\D+\Omega_{\Omega^0(E)}+
\cdots +\Omega_{\Omega^0(E)}\leqno(3.18)$$
on $\A$ which induces K\"ahler forms on $\H^\ast$ and $\Stau$.  We
will denote all these forms also by $\Omega$.  Observe that the real
gauge group $\G$ acts on $\H^\ast$ preserving $\Omega$.  As in [B-D],
Proposition 4.1, we find

\proclaim Proposition 3.19. The map $\Psi_\tau:\H^\ast\to
\lie\G$ defined by
$$\Psi_\tau\kpair=\ast\curv-{\sqrt{-1}\over
2}\sum_{i=1}^k\phi_i\otimes\phi_i^\ast+{\sqrt{-1}\over 2}\tau{\bf
I}\;,$$
is an Ad-invariant moment map for the action of $\G$ on the
symplectic manifold $(\H^\ast,\Omega)$.  Here, $\lie\G$ denotes the
Lie algebra of $\G$ and is identified with its dual via the
$L^2$-inner product.

By performing the standard infinite dimensional version of the
Marsden-Weinstein reduction (cf. [B-D], Theorem 4.5) we obtain

\proclaim Theorem 3.20.  For all values of $\tau$ satisfying
Assumption 3.12,
$\B_\tau=\B_\tau\data$ is a  K\"ahler manifold of dimension
$kd-r(k-r)(g-1)$.  Moreover, if $\tau$ is generic in the sense of
Definition 3.13, then $\B_\tau$ is compact, and is in fact a
non-singular projective variety.

\noindent  The last statement in the theorem above is a simple
generalization of the argument in [B-D-W], Theorem 6.3.
We refer to $\B_\tau\data$ as the {\it moduli space of $\tau$-stable
$k$-pairs}.

As an example, let us specialize for the moment to the case
$r=1$.  Then $\H^\ast=\H$ (Definition 3.6) is the entire space of
holomorphic pairs.  Moreover, there is no stability condition and
hence no $\tau$ dependence, once $\tau > d$.
For this case, we therefore denote
$\B(d,1,k)=\B_\tau\data$.  Moreover, comparing with \S 2   we have

\proclaim Theorem 3.21.  $\B(d,1,k)=\pspace(d,k)$ as complex
manifolds.  In particular, $\B(d,1,k)$ is a projective variety.

\noindent {\it Proof.}  To prove this, let
$pr_1 : C\times\D \lra C$
be projection onto the first factor.  Then on $\tilde U
= pr_1^\ast(E)$ there
is a tautological complex structure which is trivial in the
direction $\D$ and isomorphic to $E^\dbare$ on the slice
$C\times\{\dbare\}$.  The action of $\Gc$
lifts to $\tilde U$.  We would like to take the quotient of $\tilde
U$ by this action in order to obtain a universal bundle on $C\times
J_d$.   Unfortunately, the action of $\Gc$ on $\D$ is not free and
the constants $\C^\ast$ act non-trivially on $\tilde U$.  In rank
one, however, the choice of a point $p\in C$ allows us to express
the gauge group as a direct product
$$\Gc \simeq \Gcp \times \C^\ast\;,\leqno(3.22)$$
where
$$\Gcp =\{ g\in\Gc : g(p)=1 \}\;.$$
Then the quotient of $\tilde U$ by $\Gcp$ defines a universal
bundle $U\to C\times J_d$.  Note that the normalization (2.3) is
satisfied by this choice.  Let
$${\eqalign{
\tilde\rho : C\times\D \lra \D \cr
\rho : C\times J_d \lra J_d \cr
}}\leqno(3.23)$$
denote the projection maps.
Since we assume $d > 2g-2$, the direct images $\tilde\rho_\ast\tilde
U$ and $\rho_\ast U$ are vector bundles on $\D$ and $J_d$,
respectively.  As in (2.4), we consider
$${\eqalign{
\uktilde &=\tilde\rho_\ast\tilde
U\oplus\cdots\oplus\tilde\rho_\ast\tilde U\cr
\uk &=\rho_\ast
U\oplus\cdots\oplus\rho_\ast U\cr
}}\leqno(3.24)$$
where the fiberwise direct sums of vector bundles are taken $k$ times.
Since $\tilde\rho$ is $\Gcp$
equivariant, the quotient of the vector bundle
$\tilde U_\ast\to\D$ by $\Gcp$ is isomorphic to the
bundle $U_\ast\to J_d$, hence the quotient of the projective
bundle $\P(\uktilde)\to \D$ by $\Gcp$ is
isomorphic to $\P(\uk)=\pspace(d,k)$.
But the quotient of
 $\P(\uktilde)\to \D$ by $\Gcp$
is the same as the quotient of
$\uktilde-\{0\}$ by $\Gc$, which by
definition is the space $\B(d,1,k)$.

Now assume $r\geq 2$.  In the case $d/r < \tau <\mu_+$, where
$\mu_+$ is the smallest rational number greater than $d/r$ which can
appear as the slope of a subbundle of $E$, it is easy to see that if
$\kpair\in\Stau$, then $E^\dbare$ is semistable (cf. [B-D],
Proposition 1.7).  In this case we have the following

\proclaim Proposition 3.25.  For $r\geq 2$, $g\geq 2$ and $d/r <
\tau < \mu_+$, the natural map
$$\pi : \B_\tau\data\lra \N(d,r)$$
is a morphism of algebraic varieties,
where $\N(d,r)$ denotes the Seshadri
compactification of the moduli space of rank $r$ stable bundles of
degree $d$.

\noindent {\it Proof.}  The proof follows along the lines of [B-D],
Theorem 6.4, which uses the convergence of the gradient flow of the
Yang-Mills functional in [D].

In the case where $d$ and $r$ are coprime, there exists a universal
bundle $V\to C\times\N(d,r)$
such that $V$ restricted to $C\times\{\dbare\}$
is a stable bundle of degree $d$ isomorphic to $E^\dbare$.
Let $\rho: C\times\N\to \N$ be the
projection map.  For $d > 2r(g-1)$, the range we are considering,
the push-forward $\rho_\ast V$ is
a vector bundle on $\N$.  The map $\pi$ in Proposition
3.25 suggests the following analogue of Theorem 3.21:

\proclaim Theorem 3.26. For $r\geq 2$, $d/r <\tau <\mu_+$ and
$d$, $r$ coprime, $\B_\tau\data\simeq
\P(\vk)$ as projective varieties.

\noindent {\it Proof.}  The proof is similar to that of Theorem
3.21, only now the lack of a decomposition (3.22) makes the
construction of a universal bundle more delicate.  Indeed, if $d$,
$r$ are not coprime, such a bundle does not exist.  Let $\tilde
U\to C\times\D_s$ be the holomorphic bundle
defined as before, where now $\D_s$ denotes
the stable holomorphic structures on $E\to C$.  Since $\Gc$ does not
act freely on $\D_s$, $\tilde U$ does not descend.   However,
according to [A-B], pp. 579-580, when $d$ and $r$ are coprime we can
find a line bundle $\tilde\L\to\D_s$
with an action of $\Gc$ on which $\C^\ast\subset\Gc$
acts by multiplication on the fiber.
Lifting $\tilde\L$ to $C\times\D_s$, we define
$$\tilde V=\tilde U\otimes \tilde\L^\ast\;.\leqno(3.27)$$
Since $\C^\ast$ acts trivially on $\tilde V$, we have an action by
$\Gcbar=\Gc/\C^\ast$, and therefore $\tilde V$ descends to a
universal bundle $V\to C\times \N$.  Let $\tilde\rho$ be the lift
of the projection map $\rho$ (cf. (3.23)). Then clearly
$$ \tilde\rho_\ast\tilde V=
\tilde\rho_\ast\tilde U\otimes \tilde\L^\ast\;,\leqno(3.28)$$
so $\P(\vktilde)\simeq\P(\uktilde)$
as projective bundles on $\D_s$.  As in the proof of Theorem
3.21, the quotient of
$\P(\vktilde)$
by $\Gcbar$ gives $\P(\vk)\to\N(d,r)$, whereas the quotient of
$\P(\uktilde)$
is by definition the space $\B_\tau\data$.  This completes the
proof.

Next, we would like to show  that the spaces
$\B_\tau$ are {\it fine moduli spaces} parameterizing $\tau$-stable
$k$-pairs.  What is needed is a construction of universal bundles
$U_\tau$ on $C\times\B_\tau$ and $k$ universal sections, i.e. a map
of sheaves $\O^k_{C\times\B_\tau}\to U_\tau$.  Let
$$pr_1 : C\times\D\times\Omega^0(E)\times\cdots\times
\Omega^0(E)\lra C$$
denote projection onto the first factor.  As before,
on $pr_1^\ast(E)$
there is a tautological complex structure which is trivial in the
direction $\D\times\Omega^0(E)\times\cdots\times\Omega^0(E)$ and
isomorphic to $E^\dbare$ on the slice $C\times\kpair$. Let $\tilde
U_\tau$ denote the restriction $pr_1^\ast(E)$ to $C\times\Stau$.
There are $k$ tautological holomorphic sections
$\tilde\Phi_1,\ldots,\tilde\Phi_k$ of $\tilde U_\tau$ defined by the
property that the restriction of $\tilde\Phi_i$ to $C\times\{\dbare,
\phi_1,\ldots,\phi_k\}$ is $\phi_i$.  Next, observe that the complex
gauge group $\Gc$ acts {\it freely} on $\Stau$ and  $U_\tau$, and the
universal sections $\tilde\Phi_i$ are $\Gc$-equivariant with respect
to this action.  This implies that $\tilde U_\tau$ and
the $\tilde\Phi_i$'s descend to a bundle $U_\tau\to C\times\B_\tau$
and universal sections $\Phi_1,\ldots,\Phi_k$.  We will
denote all this by
$$\vec\Phi : \O^k_{C\times\B_\tau}\lra U_\tau\;,\leqno(3.29)$$
as mentioned above.
To summarize, we have

\proclaim Proposition 3.30.  There exists a {\it universal k-pair}
$(U_\tau\data , \vec\Phi)$ on $C\times \B_\tau\data$,
i.e. a universal rank $r$ bundle $U_\tau\data\to C\times
\B_\tau\data$ with $k$ holomorphic sections $\Phi_1,\ldots,\Phi_k$.

This universal $k$-pair, and especially the way it depends on the
parameter $\tau$, will be of fundamental importance to the
calculations in \S 5. But before closing this subsection, it will be
important to have some compatibility between the universal pair
given in Proposition 3.30 and the universal bundle $V
\to C\times\N(d,r)$ in the case where $d$ and $r$ are coprime and
$d/r < \tau < \mu_+$.  To do this, we first give an explicit
description of the anti-tautological line bundle $\O(1)$ on
$\B_\tau\data$ coming from the identification (3.26).
Let
$\tilde U$, $\tilde \L$, and $\tilde V$ be as above.  Consider the
map $\tilde \pi : \uktilde\to \D_s$.
Then the lift $\tilde\pi^\ast\tilde\L$ has a $\Gc$ action, and
therefore the quotient defines a line bundle $\L\to\B_\tau$.

\proclaim Proposition 3.31.  The bundle $\L\to\B_\tau$ defined above
is isomorphic to the anti-tauto-\break
logical line bundle $\O(1)$ under the
identification (3.26).

\noindent Note that if $\tilde\L$ is changed by a line bundle
${\cal F}\to\N$, then $V\mapsto V\otimes{\cal F}^\ast$, and hence
$\O(1)\mapsto\O(1)\otimes\pi^\ast{\cal F}$.  The proposition is thus
consistent with this fact.

\noindent {\it Proof.} It suffices
to check that the direct image $\pi_\ast\L\simeq
(\vk)^\ast$.  Since $\tilde \pi$ is $\Gcbar$ equivariant, it
suffices to check this on $\D_s$;  that is, if
$\tilde\pi^\ast\tilde\L/\C^\ast$ denotes the quotient of
$\tilde\pi^\ast\L$ on $\tilde\pi : \P(\uktilde)\to\D_s$
(here we have used $\tilde\pi$ to also denote the induced map), we
must show that
$$\tilde\pi_\ast\left(
\tilde\pi^\ast\tilde\L/\C^\ast\right)\simeq(
\vktilde)^\ast$$
as $\Gcbar$ bundles.
Suppose that $\uktilde$ has rank $N$.  Then with
respect to a local trivialization,
local
sections of $\tilde\pi^\ast\tilde\L\to\uktilde$
are functions on $
\C^N$ with values in the sheaf of local sections of $\tilde\L$.
Requiring $\C^\ast$ equivariance implies that these maps are linear.
Finally, pushing forward by $\tilde\pi$ we get an isomorphism of
sheaves
$$\tilde\pi_\ast\left(\tilde\pi^\ast\tilde\L/\C^\ast\right)\simeq(
\uktilde)^\ast\otimes \tilde\L\;.$$
The result follows from (3.28).

Now consider the universal bundle $V\to C\times\N(d,r)$ defined by
$\tilde\L$.
We have the diagram:
$$\matrix{ f^\ast V &\lra &
C\times \vk\cr
&&\mapdown{f}\cr
V &\lra &C\times\N\cr}$$
where $f$
is the identity on the first factor and the bundle projection
on the second.
The action of $\C^\ast$ on $\vk$ lifts to $f^\ast V$.  Then
it follows exactly as in the proof of Proposition 3.31 that the
quotient bundle $f^\ast V/\C^\ast$ on $C\times\P(\vk)$ is
isomorphic to $\pi^\ast V\otimes\O(1)$, where $\pi$ denotes the map
$C\times\P(\vk)\to C\times\N$.  Moreover, observe that the
tautological sections of $f^\ast V\to C\times\vk$  are
invariant under the action of $\C^\ast$, and so we obtain
universal sections
$$\vec\Psi : \O_{C\times\P(\vk)}^k
\lra\pi^\ast V\otimes\O(1)\;.$$

\proclaim Corollary 3.32.  In the case where $d/r < \tau <\mu_+$
and $d$ and $r$ are coprime, the identification (3.26) gives an
isomorphism of the $k$-pairs $(U_\tau\data,\vec\Phi)$ and $(
\pi^\ast V\otimes\O(1), \vec\Psi)$.

\noindent{\it Proof.} The isomorphism of bundles follows from the
definition of $U_\tau$, (3.27), and Proposition 3.31.  The fact
that the sections pull back is straightforward to verify.

\bigskip
\centerline{\it \S 3.3 Masterspace and flips}

\noindent In this subsection we specialize to the case $r=2$ and
examine the dependence of the moduli spaces $\B_\tau$ on the
parameter $\tau$ within the admissible range $(d/2, d-(2g-2))$.  More
precisely, we show as in [B-D-W] that when we pass an integer value in
$(d/2,d-(2g-2))$, i.e. the non-generic values in the sense of
Definition 3.13, the spaces $\B_\tau$ are related by a {\it flip},
by which we mean simply a blow-up followed by a blow-down of the
exceptional divisor in a different direction.  We show this by
directly generalizing the ``master space" construction developed in
[B-D-W].  Since this part is quite straightforward we shall be brief
and refer to the latter paper for more details.
Finally, we describe the relationship between the
universal bundles $U_\tau$ on the various $\B_\tau$'s. This will be
important for the calculations in \S 5.3.

Recall that the space $\H^\ast$ of rank two
holomorphic $k$-pairs
$\kpair$ satisfying $\mu_M(E) < d-(2g-2)$ has the structure of a
K\"ahler manifold and is acted on holomorphically and symplectically
by the gauge group $\G$.  We may define a character on $\G$ as
follows (cf. [B-D-W], \S 2): First, choose a splitting
$\G\simeq\G_1\times\Upsilon$, where $\Upsilon$ is the group of
components of $\G$ and $\G_1$ denotes the connected component of the
identity.  For $g\in\G_1$ we define $\chi_1(g)$ to be the unique
element of U(1) such that $\det g=\chi_1(g)\exp u$, where $u:C\to
i{\Bbb R}$ satisfies $\int_C u=0$.  We extend to a character $\chi$
on $\G$ by letting $\chi_1$ act trivially on $\Upsilon$.  Let $\G_0$
denote the kernel of $\chi$, i.e. we have an exact sequence of
groups
$$1\lra\G_0\lra\G\mapright{\chi}\hbox{U(1)}\lra 1\;.\leqno(3.33)$$
The new space $\hat\B$ is then obtained by symplectic reduction of
$\H^\ast$ by the action of the smaller group $\Go$.  Specifically,
by Proposition 3.19,
$$\Psi\kpair =\ast\curv -{\sqrt{-1}\over
2}\sum_{i=1}^k\phi_i\otimes\phi_i^\ast\leqno(3.34)$$
is a moment map for the action of $\G$ on $\H^\ast$. Hence, as in
[B-D-W], Proposition 2.7,
$$\Psi_0\kpair=\Psi\kpair-{1\over 2}\int_C\tr\Psi\kpair\cdot{\bf
I}\leqno(3.35)$$
is a moment map for the action of $\G_0$.
Since $r=2$, one can show that
$\G_0$ and its complexification $\Gc_0$ act freely on
$\H^\ast$ (cf. [B-D-W], Proposition 2.19).
It then follows that
$$\hat\B=\Psi_0^{-1}(0)/\G_p\;,\leqno(3.36)$$
is a smooth, symplectic manifold.
A small variation of Proposition 3.20 along
the lines of [B-D-W], Theorem 2.16
then proves that $\hat\B$ is a complex
manifold  with the symplectic structure giving $\hat\B$ a K\" ahler
manifold structure.
Finally, observe that there is a residual $\G/\G_0\simeq\hbox{ U(1)}$
action
on $\hat\B$ that is quasifree, i.e. the stabilizer of every point is
either trivial or the whole U(1).  Moreover, this action is clearly
holomorphic and symplectic.  Let
$f:\hat\B\to{\Bbb R}$
denote the associated moment map.
As in [B-D-W], it follows that $f$ is given by
$$f\kpair={1\over 8\pi}\sum_{i=1}^k\Vert\phi_i\Vert^2_{L^2}+{d\over
2}\;.\leqno(3.37)$$
The generic (i.e. nonintegral) values of $\tau \in (d/2,d-(2g-2))$
correspond to the level sets where the U(1) action is free, hence
the reduced spaces $f^{-1}(\tau)/\hbox{U(1)}$ for such $\tau$ inherit a
K\"ahler manifold structure, and indeed
$$\B_\tau\simeq f^{-1}(\tau)/\hbox{U(1)}$$
as K\"ahler manifolds.  We summarize the above discussion by the
following

\proclaim Theorem 3.38. There is a K\"ahler manifold $\hat\B$ with a
holomorphic, symplectic, quasifree U(1) action whose associated
moment map is given by (3.37).  Moreover, with the induced K\"ahler
structure from $\hat\B$, $f^{-1}(\tau)/\hbox{U(1)}
\simeq\B_\tau$ as K\"ahler
manifolds for any noninteger value $\tau\in(d/2,d-(2g-2))$.

The space $\hat\B$ is key to understanding the relationship between
the $\B_\tau$'s for different values of $\tau$.  As explained in
[B-D-W], this picture may best be understood via the Morse theory of
the function $f$ on $\hat\B$.  First, observe that the critical
values of $f$ are precisely the non-generic values of $\tau$, i.e.
the integers in $(d/2, d-(2g-2))$.  For generic
$\tau <\tau'$, if $[\tau,\tau ']\cap {\Bbb Z}=\emptyset$,
then the gradient flow of $f$ gives a biholomorphism $\B_\tau\simeq
\B_{\tau'}$, and indeed by Definition 3.3, $\B_\tau
=\B_{\tau'}$ as complex manifolds (though not as K\"ahler
manifolds). Next, suppose that $\tau$ is a critical value of $f$,
and let $Z_\tau$ denote the critical set corresponding to $\tau$.
As in [B-D-W], Example 3.5, It follows from Theorem 3.5  that
$$Z_\tau\simeq \B(d-\tau, 1,k)\times J_\tau \;,\leqno(3.39)$$
where as in \S 2, $J_\tau$ denotes the degree $\tau$
Jacobian variety of $C$.   Let $W^+_\tau$ ($W^-_\tau$) denote the
stable (unstable) manifolds of gradient flow by the function $f$.
We can express $W^\pm_\tau$ in terms of
extensions of line bundles as follows (cf. [B-D-W], \S 4):

\proclaim Proposition 3.40.  (i) $W^+_\tau$ consists of stable
$k$-pairs $(E,\vec\phi)$ such that $E$
fits into an exact sequence
$$0\mapright{} F\mapright{} E\mapright{\pi} Q_\phi\mapright{}0$$
such that under the isomorphism (3.39), $(Q_\phi,\pi(\vec\phi))\oplus
F\in Z_\tau$, and $\tau$ is maximal with respect to this property.
(ii)  $W^-_\tau$ consists of stable $k$-pairs $(E,\vec\phi)$ such
that $E$ fits into an exact sequence
$$0\mapright{} E_\phi\mapright{} E\mapright{} F\mapright{} 0$$
such that under the isomorphism (3.39), $(E_\phi,\vec\phi)\oplus
F\in Z_\tau$, and $\tau$ is minimal  with  respect to this property.

\noindent {\it Proof.}  Immediate generalization of [B-D-W],
Propositions 4.2 and 4.3.

Pick $\tau\in (d/2, d-(2g-2))\cap{\Bbb Z}$ and $\varepsilon > 0$
sufficiently small such that $\tau$ is the only integer value in
$[\tau-\varepsilon,\tau+\varepsilon]$.  Let
$\P(W^\pm_\tau)$ be the images of $W^\pm_\tau\cap
f^{-1}(\tau\pm\varepsilon)$ under the quotient map
$$f^{-1}(\tau\pm\varepsilon)\lra
\B_{\tau\pm\varepsilon}=f^{-1}(\tau\pm\varepsilon)/\hbox{U(1)}\;.$$
It follows from Proposition 3.40 that
$$\sigma_\pm : \P(W_\tau^\pm)\lra Z_\tau\leqno(3.41)$$
are projective bundles over the critical set.  Furthermore, by
direct application of the Morse theory and the description
of Guillemin and Sternberg [G-S] one can show as in [B-D-W],
Theorem 6.6,

\proclaim Theorem 3.42.  Let $\tau$ and $\varepsilon$ be as above.
Then there is a smooth
projective variety $\tilde \B_\tau$ and holomorphic
maps $p_\pm$
$$\matrix{&& \tilde\B_\tau &&\cr
&{}^{p_-}\swarrow &&\searrow {}^{p_+}&\cr
\B_{\tau-\varepsilon}&&&& \B_{\tau+\varepsilon}\cr}$$
Moreover, $p_\pm$ are blow-down maps onto the smooth subvarieties
$\P(W^\pm_\tau)$,
and the exceptional divisor $A\subset
\tilde \B_\tau$ is the fiber product
$$ A\simeq \P(W^-_\tau)\times_{Z_\tau}
\P(W^+_\tau)\;.$$

We end this section with some important facts concerning the
universal bundles $U_\tau$ and the normal bundles to the centers
$\P(W^\pm_\tau)$.

\proclaim Proposition 3.43.  The normal bundle
$\nu\left(
\P(W^\pm_\tau)\right)$ of
$\P(W^\pm_\tau)$  in $\B_{\tau\pm\varepsilon}$ is given
by
$$\nu\left(
\P(W^\pm_\tau)\right)
=\sigma_\pm^\ast W_\tau^\mp \otimes \O_{\P(W^\pm_\tau)}(-1)\;,$$
where $\sigma_\pm : \P(W^\pm_\tau)
\to Z_\tau$ is the projective bundle (3.41) associated to the stable
and unstable manifolds, and $\O_{\P(W^\pm_\tau)}
(-1)$ are the tautological line bundles on
$\P(W^\pm_\tau)$.

\noindent {\it Proof.} Let $\hat\sigma_\pm$ denote the maps
$$\hat\sigma_\pm : W_\tau^\pm\cap f^{-1}(\tau\pm\varepsilon)\lra
Z_\tau$$
induced by the gradient flow of $f$.  Since the flow is invariant
under the circle action, $\hat\sigma_\pm$ lift $\sigma_\pm$.
We focus on $W^+$, the argument for $W^-$ being exactly the same.
Consider the diagram
$$\matrix{\hbox{U(1)}&&\cr
	\mapdown{}&&\cr
	W_\tau^+\cap f^{-1}(\tau+\varepsilon)
	&\mapright{\hat\sigma_+}&Z_\tau
	\cr
	\mapdown{}&&\mapdown{}\cr
	\P(W_\tau^+)&\mapright{\sigma_+}&Z_\tau\cr}$$
Since $\hat\sigma_+$ is a retract, the normal bundle to
$W^+_\tau\cap f^{-1}(\tau+\varepsilon)$ in
$f^{-1}(\tau+\varepsilon)$ is the pullback by $\hat\sigma_+$
of a bundle on $Z_\tau$.  Since the tangent bundle to $\hat\B$ along
$Z_\tau$ decomposes under the U(1) action into $TZ_\tau\oplus
W^+\oplus W^-$, the bundle in question is clearly $W^-$.  Now an
element $e^{i\theta}\in\hbox{U(1)}$ acts on $W^-$ by $e^{-2i\theta}$
and on $W^+$ by $e^{2i\theta}$.  Hence upon taking quotients,
$\hat\sigma_+^\ast W^-/\hbox{U(1)}$ is just $\sigma_+^\ast
W^-$ twisted by the tautological line bundle.  This proves the
topological equivalence of the bundles.  The holomorphic equivalence
follows, since the U(1) action and the gradient flow are both
holomorphic.

Let $\L_\phi$ and $\L_s$ denote the pullbacks to $C\times Z_\tau$ of
the universal bundles on $C\times\B(d-\tau,1,k)$ and $C\times
J_\tau$ under the identification (3.39).

\proclaim Proposition 3.44.  The following are exact sequences of
holomorphic bundles on $C\times\P(W_\tau^\pm)$:
$$0\lra\sigma_-^\ast\L_\phi\lra
U_{\tau-\varepsilon}\bigr|_{C\times\P(W_\tau^-)}\lra\sigma_-^\ast\L_s
\otimes\O_{\P(W_\tau^-)}(-1)\lra 0\;;\leqno(i)$$
$$0\lra\sigma_+^\ast\L_s\otimes\O_{\P(W_\tau^+)}(+1)\lra
U_{\tau+\varepsilon}\bigr|_{C\times\P(W_\tau^+)}\lra\sigma_+^\ast
\L_\phi \lra 0\;.\leqno(ii)$$

\noindent{\it Proof.}  Recall from the construction Proposition 3.30
of the universal bundles  that $U_\tau$ is the quotient by the
action of $\Gc$  of the restriction of the natural bundle on
$C\times\H^\ast$.  By quotienting out by the smaller group $\Gc_0$
we obtain a universal bundle $\hat U$ on $C\times\hat\B$ which
descends under the U(1) reduction at $\tau$ to $U_\tau$.  With this
understood, consider part (i) of the proposition.  By the above
construction, the universal bundle $\hat U$  restricted to $Z_\tau$
canonically splits
$$\hat U\bigr|_{C\times Z_\tau}\simeq \L_\phi\oplus\L_s\;.$$
Again since $\hat\sigma_-$ is a retract, we have a natural sequence
$$0\lra\hat\sigma_-^\ast\L_\phi\lra
\hat U\bigr|_{C\times f^{-1} (\tau-\varepsilon)}\lra\hat
\sigma_-^\ast\L_s \lra 0\;.\leqno(3.45)$$
Moreover, (3.45) is an exact sequence of ``CR-bundles," by which we
mean that the $\dbar_b$-operator on $\hat\sigma_-^\ast\L_\phi$ is
induced from the $\dbar_b$-operator on $\hat U\bigr|_{C\times
f^{-1}(\tau-\varepsilon)}$, and likewise for
$\hat\sigma_-^\ast\L_s$.  Note that this statement would not be true
if we reversed the sequence.  One can determine the quotient of
(3.45) by the U(1) action if the action is trivialized on the
$\L_\phi$ part.  Specifically, the action by $e^{i\theta}$ is gauge
equivalent (in $\G_0$) to the action by $g_\theta=\hbox{diag}(1,
e^{2i\theta})$.  Now $g_\theta$ acts trivially on
$\hat\sigma_-^\ast\L_\phi$ and by $e^{2i\theta}$ on
$\hat\sigma_-^\ast\L_s$.  Since $g_\theta$ also acts as
$e^{-2i\theta}$ on $W^-$, we see that
$$\hat\sigma_-^\ast\L_s/\hbox{U(1)}\simeq
\sigma_-^\ast\L_s\otimes\O_{\P(W_\tau^-)}(-1)\;.$$
As before, this proves part (i) as a topological statement, and the
holomorphicity follows from the holomorphicity of the U(1) action.
Part (ii) is proved similarly.

\beginsection 4.  Maps to Grassmannians

\bigskip
\centerline{\it \S 4.1 The Uhlenbeck compactification}

\noindent The main purpose of this section is to show that the Uhlenbeck
compactification of $\M\data$, the space of holomorphic maps of
degree $d$ from $C$ to $G(r,k)$, admits the structure of a projective
variety.  In the case of maps into projective space, it is actually
a smooth, holomorphic projective bundle over the Jacobian.

We first introduce $\U\data$ set theoretically.  Let $0\leq l\leq
d$ be an integer, and let $C_l$ denote the $l$-th symmetric product
of the curve $C$.  We set
$$\U\data=\bigsqcup_{0\leq l\leq d}C_l\times\M(d-l,r,k)
\;,\leqno(4.1)$$
and we will denote elements of $\U\data$
by ordered pairs $(D,f)$.  Given $(D,f)$ we associate the
distribution $e(f)+\delta_D$, where $e(f)=|df|^2$ is the energy
density of the map $f$ with respect to the fixed K\"ahler metrics on
$C$ and $G(r,k)$, and $\delta_D$ is the Dirac distribution supported
on $D$.  We topologize $\U\data$
by giving a local basis of neighborhoods
around each point $(D,f)$ as follows:  Pick a basis of neighborhoods
$N$ of $f\in\M(d-l,r,k)$ in the $C^\infty_0(C\setminus D)$ topology
and a basis of neighborhoods $W$ of $e(f)+\delta_D$ in the
weak*-topology.  Set
$$V(D,f)=\left\{ (D',f')\in\U\data : f'\in N\hbox{ and
}e(f')+\delta_{D'}\in W\right\}\;.\leqno(4.2)$$
Since both the $C^\infty_0$ and weak* topologies are first
countable, so is the topology defined on $\U\data$.  In terms of
sequences then, $(D_i,f_i)\to (D,f)$ in $\U\data$ if and only if:
\itemitem{(i)} $f_i\to f$ in the $C^\infty_0(C\setminus D)$
topology, and
\itemitem{(ii)} $e(f_i)+\delta_{D_i}\to e(f)+\delta_D$ in the
weak*-topology.

\noindent  By the theorem of Sacks and Uhlenbeck [S-U], $\U\data$
with this topology is compact.

\proclaim Definition 4.3.  The space $\U\data$ with topology
described above is called the {\it Uhlenbeck compactification} of
$\M\data$.

Let us consider for the moment the case $r=1$.
Then the Uhlenbeck compactification of holomorphic maps of degree
$d$ from $C$ to $\P^{k-1}$ is defined once we have chosen K\"ahler
metrics on $C$ and $\P^{k-1}$. For $C$ we use the metric of \S 3,
and for $\P^{k-1}$ we choose the Fubini-Study metric.
Then we have

\proclaim Theorem 4.4.  If $r=1$,
then $\U(d,1,k)$ is homeomorphic to $\B(d,1,k)$.  In
particular, the Uhlenbeck compactification for holomorphic maps of
degree $d> 2g-2$ into projective space has the structure of a {\it
non-singular} projective variety.

\noindent The proof of Theorem 4.4 is somewhat lengthy,
so we will split it into several lemmas.
Our first goal is to define a map
$$ u : \B(d,1,k)\lra \U(d,1,k)\;.\leqno(4.5)$$
This is achieved by the following

\proclaim Lemma 4.6.   Given a stable $k$-pair $(L,\vec
\phi)\in\B(d,1,k)$
there is a unique divisor $D\in C_l$ and holomorphic map
$$f_{(L,\vec\phi)} : C\setminus D\lra \P^{k-1}$$
which extends uniquely to a holomorphic map of degree $d-l$ on $C$.

\noindent {\it Proof.}  Given $(L,\vec\phi)\in\B(d,1,k)$, recall the
map $C\to\P^{k-1}$ given by (2.6).
This defines an algebraic map on $C\setminus D$,
where $D=\{ p\in C : \vec\phi(p)=0\}$,
counted with multiplicity.
By the properness of $\P^{k-1}$ this map extends to a holomorphic
map
$$f_{(L,\vec\phi)}: C\longrightarrow \P^{k-1}\leqno(4.7)$$
of degree $d-l$ where $l=\deg D$.

\proclaim Definition 4.8.  The map $u$ of (4.5) is defined by
setting $u[L,\vec\phi]=(D, f_{(L,\vec\phi)})$, where $D$ and $f_{(L,
\vec\phi)}$ are defined by Lemma 4.6.

In order to prove Theorem 4.4 it is convenient first to write a local
expression for the map (4.7) in terms of the
homogeneous coordinates of $\P^{k-1}$.
Let $z$ be a local coordinate centered at $p\in D$, and suppose that
$m$ is the minimal order of vanishing of the $\phi_i$'s at $p$.
Then in a deleted neighborhood of $p$, the map (2.6) is clearly
equivalent to
$$f_{(L,\vec\phi)}(z)=\left[{\phi_1(z)\over z^m},\ldots,{\phi_k(z)\over
z^m}\right]\;.\leqno(4.9)$$
Since for some $i$, $\lim_{z\to 0}\phi_i(z)/z^m\neq 0$, (4.9) is the
desired extension.

After this small digression, we continue with

\proclaim Lemma 4.10.  The map $u$ is a bijection.

\noindent {\it Proof.} A set theoretic inverse to $u$ can be
constructed as follows:  Let $S^\ast\to\P^{k-1}$ denote the
anti-tautological bundle with $k$ tautological sections $z_1,\ldots,
z_k$ given by the coordinates of $\C^k$.  Given $(D,f)\in\U(d,1,k)$,
set $L=f^\ast S^\ast\otimes\O_C(D)$, where $\O_C$ denotes the structure
sheaf of $C$, and $\phi_i=f^\ast z_i\otimes 1_D$, $i=1,\ldots, k$,
where $1_D$ denotes a choice of holomorphic section of $\O_C(D)$
vanishing at precisely $D$.  Then it is clear that
$u^{-1}(D,f)=[L,\phi_1,\ldots,\phi_k]$.  Observe that a different
choice of $1_D$ amounts to a rescaling of the $\phi_i$'s and so does
affect the definition of $u^{-1}$.

The next step is to prove

\proclaim Lemma 4.11.  The map $u$ is continuous.

\noindent {\it Proof.}  Assume $[\dbar_n,\vec\phi_n]\to
[\dbar,\vec\phi]$ in $\B(d,1,k)$.  Write $(D_n,
f_n)=u[\dbar,\vec\phi_n]$ and
$(D,f)=u[\dbar,\vec\phi]$.  Since
$\pi:\B(d,1,k)\mapright{\pi} J_d$ is continuous, we may assume that
the operators $\dbar_n\to\bar\partial$ on $L$ and that we can
choose  local holomorphic trivializations simultaneously for
all $\dbar_n$, $\dbar$.  According to the definition of the
Uhlenbeck topology we must show
\itemitem{(i)} $f_n\to f$ in $C^\infty_0(C\setminus D)$, and
\itemitem{(ii)} $e(f_n)+\delta_{D_n}\to e(f)+\delta_D$ as
distributions.

\noindent For (i), pick $p\in C\setminus D$.  Clearly, $p\in
C\setminus D_n$ for $n$ sufficiently large.  Then by
the construction in Lemma 4.6,
$$\eqalign{ f_n(p)&=[\phi_{1,n}(p),\ldots,\phi_{k,n}(p)]\cr
 f(p)&=[\phi_1(p),\ldots,\phi_k(p)]\;,\cr}$$
 where $\phi_{i,n}\to \phi_i$ smoothly as $n\to \infty$ for all
 $i=1,\ldots, k$.  This clearly implies (i).

To prove part (ii), recall that the Fubini-Study metric
on $\P^{k-1}$ is given by
$$\omega ={\sqrt{-1}\over 2\pi}\partial\dbar\log
\left(\sum_{i=1}^k|z_i|^2\right)\;.$$
Since $f_n$ and $f$ are holomorphic,
$$\eqalign{
e(f_n)(p)&=f_n^\ast\omega= {\sqrt{-1}\over 2\pi}\partial\dbar\log
\left(\sum_{i=1}^k|\phi_{i,n}(p)|^2\right)\cr
e(f)(p)&=f^\ast\omega= {\sqrt{-1}\over 2\pi}\partial\dbar\log
\left(\sum_{i=1}^k|\phi_i(p)|^2\right)\;.\cr}$$
Let $z$ be a local coordinate about a point $p\in D$, and let $V$ be
a neighborhood of $p$ satisfying $V\cap D=\{p\}$.  Furthermore, by
induction on $d$ we assume $p\not\in D_n$ for $n$ large and that the
multiplicity $m$ of $p$ in $D$ is $m\geq 1$.
Let $g\in C^\infty_0(V)$.  Then
$$\eqalign{
\int_V g e(f_n) &= \int g {\sqrt{-1}\over 2\pi}\partial\dbar\log
\left(\sum_{i=1}^k|\phi_{i,n}(z)|^2\right)\cr
&=
\int g {\sqrt{-1}\over 2\pi}\partial\dbar\log
|z|^{2m}+
\int g {\sqrt{-1}\over 2\pi}\partial\dbar\log
\left(\sum_{i=1}^k\left|{\phi_{i,n}(z)\over z^m}\right|^2\right)\cr
&=m\delta_p(g)+
{\sqrt{-1}\over 2\pi}\int_V \partial\dbar g \log
\left(\sum_{i=1}^k\left|{\phi_{i,n}(z)\over z^m}\right|^2\right)\cr}
\leqno(4.12)$$
On the other hand, since $\phi_{i,n}(z)\to z^m\phi_i(z)$ in
$C^\infty_0(V)$ for all $i=1,\ldots, k$, it follows from the
dominated convergence theorem and the fact that $\log|z|^{-m}\in
L^1(V)$ that
$$\log
\left(\sum_{i=1}^k\left|{\phi_{i,n}(z)\over z^m}\right|^2\right)
\lra \log
\left(\sum_{i=1}^k|\phi_i(z)|^2\right)$$
in $L^1(V)$. By taking limits in (4.12) we obtain
$$\int_V g e(f_n)\lra m\delta_p(g)+\int_V g e(f)\;.$$
Now by covering $C$ with $V$'s as above and using partitions of
unity we obtain the convergence (ii).  This completes the proof of
Lemma 4.6.

\noindent {\it Proof of Theorem 4.4.}  According to Lemmas 4.10 and
4.11, the map $u:\B(d,1,k)\to\U(d,1,k)$ is a continuous bijection of
compact topological spaces, and hence is a homeomorphism.  It
follows that $\U(d,1,k)$ inherits the projective bundle structure of
$\B(d,1,k)$.

Next,  we proceed to show that $\U\data$ has the structure of a
projective variety for any $r$.  Actually, since we are interested
in computing intersections, we will have to be more precise and
define a {\it scheme} structure on $\U\data$.  The reasons for this
will become apparent in the following subsection.

To begin, note that the Pl\"ucker embedding
$$G(r,k)\hookrightarrow
\P^{N-1}\;,\leqno(4.13)$$
where  $N={k\choose r}$,
induces an inclusion on the Uhlenbeck spaces
$$\U\data\hookrightarrow\U(d,1,N)\;.\leqno(4.14)$$
The next proposition is immediate from the definition of the
Uhlenbeck topology:

\proclaim Proposition 4.15.  The inclusion (4.14) is a homeomorphism
of $\U\data$ onto a closed subspace of $\U(d,1,N)$.

We now show that the image of  (4.14) has a natural
scheme structure.   We shall use the identification
$\U(d,1,N)\simeq\B(d,1,N)$ coming from Theorem 4.4.
Let ${\cal H}={\cal H}(d,1,N)$ denote the space
of holomorphic $k$-pairs $(L,\phi_1,\ldots, \phi_N)$, where $L$ is a
holomorphic line bundle of degree $d$.  On $C\times {\cal H}$ we
have the universal line bundle $\tilde
U$.  Suppose $\Xi$ is quadratic form
on ${\Bbb C}^N$, and denote by $Q_\Xi$ the quadric hypersurface in
${\Bbb P}^{N-1}$ defined by $\Xi$.  Then $\Xi$ determines a section
$\tilde\psi_\Xi$ of $\tilde U^{\otimes 2}$ as follows:  For a point
$x=(p;L,\vec\phi)\in C\times{\cal H}$, let
$$\tilde\psi_\Xi(x) =\Xi(\vec\phi(p),\vec\phi(p))\;.$$
Clearly, $\tilde\psi_\Xi$ is holomorphic, and since $\Xi$ is quadratic,
$\tilde\psi_\Xi$ is a section of $U^{\otimes 2}$.  Moreover, it follows
by definition of the action of ${\frak G}^{\Bbb C}$ on $\tilde U$ that
$\tilde\psi_\Xi$ is equivariant with respect to this action.  Therefore,
$\tilde\psi_\Xi$ descends to a holomorphic section of $U^{\otimes 2}\to
C\times \B(d,1,N)$.  We denote this section by $\psi_\Xi$.

Given a point $p\in C$, let $Z_\Xi(p)$ denote the zero scheme of
$\psi_\Xi\bigr|_{\{p\}\times\B(d,1,N)}$ in $\B(d,1,N)$.
Then we have the following

\proclaim Lemma 4.16.  Let $p_1,\ldots, p_m$ be distinct points in $C$
with $m\geq 2d+1$.  Let $Z_\Xi$ be the scheme theoretic intersection
$Z_\Xi(p_1)\cap\cdots\cap Z_\Xi(p_m)$.  Then $Z_\Xi$ corresponds set
theoretically to the set of $(D,f)\in \U(d,1,N)$ where
$f(C)\subset Q_\Xi$.

\noindent {\it Proof.}  If $f(C)\subset Q_\Xi$, then the point
$x=[L,\vec\phi]$ corresponding to $f$ satisfies $\vec\phi(p)=0$ if
$p\in D$ and $[\phi_1(p),\ldots,\phi_N(p)]\in\ker\Xi$ otherwise; in
particular, $x\in Z_\Xi$.  Conversely, suppose $x\in Z_\Xi$, and let
$(D,f)$ be the point in $\U(d,1,N)$ corresponding to $x$.  Now
$Z_\Xi(p_i)$ consists of all points $[L,\vec\phi]$ where either
$\vec\phi(p_i)=0$ or  $[\phi_1(p),\ldots,\phi_N(p)]\in\ker\Xi$.  If
$x\in Z_\Xi(p_i)$ then the former condition implies $p_i\in D$ and
the latter implies $f(p_i)\in Q_\Xi$.  Either way, $f$ maps at least
$m-l$ points into $Q_\Xi$, where $l$ is the degree of $D$. Since $f$
has degree $d-l$ and $m-l > 2(d-l)$, Bezout's Theorem implies
$f(C)\subset Q_\Xi$.  This completes the proof.

The embedding (4.13)
realizes $G(r,k)$ as the common zero locus of quadratic forms
$\Xi_i$, $i=1,\ldots,N$ (see [G-H], p. 211).
Therefore, Lemma  4.16
immediately implies

\proclaim Theorem 4.17.  The image of $\U\data$ in
$\B(d,1,N)$ is precisely the intersection \break
$Z_{\Xi_1}\cap\cdots \cap Z_{\Xi_N}$.
In particular, $\U\data$ has the structure of a projective
scheme.

\bigskip
\centerline{\it \S 4.2 The Grothendieck Quot scheme}

\noindent In this section, we will exhibit a different compactification
of the space $\M$  in terms of a certain Grothendieck Quot scheme.
This compactification is perhaps more natural in the algebraic
category and will be essential for our computations in \S 5.

Let $\F$ be a coherent sheaf on our fixed Riemann surface $C$.
As before, we denote the structure sheaf of $C$ by $\O_C$.
For each $t\in {\Bbb Z}$, the {\it coherent sheaf Euler
characteristic}
$$h_{\cal F}(t) := \chi (C,{\cal F}(tp)) =
h^0(C,{\cal F}(tp)) - h^1(C,{\cal F}(tp))$$
does not depend upon the choice of a point $p\in C$, and $h_{\cal
F}(t)$ is
a polynomial in $t$ (see [Gro]). This is referred to as the
Hilbert polynomial of the sheaf $\F$.
For example, if $E$ is a vector bundle of degree $d$ and rank $r$ on
$C$, then by the Riemann-Roch theorem,
$h_E(t) = d+rt-r(g-1)$. In particular, both the rank and degree of a
vector bundle are determined by its Hilbert polynomial.

Recall that on the Grassmannian $G(r,k)$, there is the tautological
exact sequence
$$0\lra S\lra \O_G^k \lra Q\lra 0\;,\leqno(4.18)$$
where $\O_G^k$ is the trivial bundle of rank
$k$ on $G(r,k)$
and $S$ and $Q$ are the universal bundles of rank $r$ and $k-r$,
respectively.
If $f:C\to G(r,k)$ is a holomorphic map of degree $d$, then the
pullback of the tautological quotient yields a quotient
$\O_C^k \to f^*Q \to 0$ of vector bundles on
$C$.  Furthermore, since $f^*S$ is of rank $r$ and
degree $-d$, the Hilbert polynomial
$$h_{f^*Q}(t) = h_d(t)  := kh_{\O_C}(t)
- [rt-(d+r(g-1))]\;.$$

Actually, to be more precise, the map $f$ determines an equivalence
class of quotients $\O_C^k \to F \to 0$, where
two such quotients are equivalent if there is an isomorphism
of the $F$'s which carries one quotient to the
other. Since such a quotient also clearly determines a
map from $C$ to $G(r,k)$, we have
the following

\proclaim  Lemma 4.19.  The set of degree $d$ holomorphic
maps $f:C\to G(r,k)$ may be identified with the set of
equivalence classes of quotients
$\O_C^k\rightarrow F \rightarrow 0$, where
$h_F(t) = h_d(t)$ is as defined above.

The idea behind the Quot scheme compactification is to expand
the set of quotients to include quotients $\O_C^k \to
\F \to 0$, where $\F$ is a coherent sheaf with
Hilbert polynomial $h_\F(t) = h_d(t)$. Following Grothendieck
(see [Gro]), one considers the contravariant
``quotient'' functor assigning to
each scheme $X$ the set of quotients $\O_{C\times X}^k
\to \widetilde \F\to 0$
(modulo equivalence) of coherent sheaves
on $C\times X$ such that $\widetilde \F$ is flat over $X$ with
relative Hilbert polynomial $h_d(t)$. Grothendieck's theorem is the
following:

\proclaim Theorem 4.20. The quotient functor is representable by a
projective scheme. That is, there is a projective scheme $\Q
=\Q\data$ together with a universal quotient $\O_{C\times
\Q}^k \to \widetilde \F \rightarrow 0$
flat over $\Q$ with relative Hilbert polynomial $h_d(t)$,
such that each of the flat quotients over $X$
defined above is equivalent to the pullback of the
universal quotient under a unique morphism from $X$ to $\Q$.

The projective scheme $\Q$ defined above is clearly uniquely determined
(up to isomorphism)
and is called the {\it Grothendieck Quot scheme}. Of course, the closed
points of $\Q$ correspond to equivalence classes of quotients on
$C$, so the scheme $\Q$ parametrizes such equivalence classes
and by Lemma 4.19, the subset of $\Q$
corresponding to vector bundle quotients
parametrizes holomorphic maps $f:C\to G(r,k)$.

But every point of $\Q$ determines a unique holomorphic
map to the Grassmannian via the following

\proclaim Lemma 4.21.  (i) The kernel of a quotient $\O^k
\to \F\to 0$ is always a vector bundle
on $C$. (ii) If $\O_{C\times X}^k
\to\widetilde \F \to 0$ is a flat quotient
over $C\times X$, then the kernel is locally free.

\noindent {\it Note:} For schemes
$X$ that are not necessarily
smooth, we will use the terminology ``locally free'' instead of ``vector
bundle''.

Part (i) of Lemma 4.21 tells us that a quotient $\O^k
\to \F \to 0$ induces an injection of sheaves
$0\to E^*\to \O^k$, where $E$ is a vector
bundle on $C$ of rank $r$ and degree $d$. Moreover, an equivalence
class of quotients clearly corresponds to the analogous
equivalence class of
injections. But an injection of sheaves induces an injection of the
fibers at all but a finite number of points, hence a rational map of
$C$ to the Grassmannian.
As in Lemma 4.6, this in turn defines a holomorphic map of lower
degree from $C$ to $G(r,k)$.

We may therefore interpret the Quot scheme as a fine moduli space
for the ``injection'' functor assigning to $X$ the set of sheaf
injections $0\to\widetilde E^*\to\O_{C\times
X}^k$, and dualizing, we get the following corollary to Theorem 4.2:

\proclaim Corollary 4.22.
The Quot scheme $\Q$ is a fine moduli
space for the functor which assigns to each scheme $X$ the set of
sheaf maps $\O_{C\times X}^k \to \widetilde E$
(modulo equivalence) subject to the following conditions:
$\widetilde E$ is locally free, for each closed point
$x\in X$ the restriction of $\widetilde E$ to
$C\times \{x\}$ has rank $r$ and
degree $d$, and the restriction of the sheaf map
to $C\times \{x\}$ is surjective at
all but a finite number of points.

If we compare Corollary 4.22 and Definition 3.3, we see that
we may interpret
the Quot scheme as a fine moduli space for $\tau$-stable
$k$-pairs if $\tau > d$
(cf. Proposition 3.14).
If $r \geq 2$,  these values of $\tau$ are outside the admissible
range, so the Quot scheme does not coincide with
one of the smooth moduli spaces constructed in \S 3.
Indeed, the Quot scheme has singularities in general.
However, in the rank one case, we have the following:

\proclaim  Corollary 4.23.
The Quot scheme $\Q(d,1,k)$ is isomorphic to the
projective bundle $\pspace(d,k)$ over $J_d$ defined in (2.4).

\noindent {\it Proof.} By Corollary 4.22 and Theorem 3.5,
the moduli space of stable $k$-pairs for rank one bundles
is a fine moduli space representing the same functor as
$\Q(d,1,k)$. Therefore they are isomorphic, and by Theorem 3.22,
$\B(d,1,k)$ is isomorphic to $\pspace(d,k)$.

Putting Corollary 4.23 together with Theorem 4.16, we get the following:

\proclaim  Theorem 4.24.
 There is an algebraic surjection $u: \Q\data
\to \U\data$ which is an isomorphism on $\M\data$.

\noindent {\it Proof.} Let $N={k\choose r}$ and
consider the universal sheaf map
from Corollary 4.21:
$\vec\Psi: \O_{C\times \Q}^k \to \widetilde E$.
The map
$$\wedge^r(\Psi): \O_{C\times \Q}^N \lra\bigwedge^r\widetilde E\;,$$
determines a morphism
$$w: \Q\data \lra \Q\left(d,1,N\right)\;.$$
Using Corollary 4.23, the image of $w$ is easily seen to be
precisely $\U\data$.
Now the morphism $w$ is not in general an embedding.
However, if $x\in \Q\data$ parametrizes a  surjective map $\O_C^k\to
E$, then as we saw above the quotient is completely determined
by the corresponding map to $G(r,k)$. Via the Pl\"ucker
embedding of the Grassmannian, the point $x$ is recovered from
$w(x)$. Thus $w$ is an embedding when restricted to $\M\data$,
which completes the proof of the theorem.

If $D$ is an effective divisor on $C$
of degree $\delta$, let $\widetilde E(D)$
denote the tensor product $\widetilde E\otimes \pi^*\O_C(D)$
of bundles on $C\times \Q$. If we let $\tilde d = d +
r\delta$, then the natural sheaf map $\O_C\to
\O_C(D)$ induces a map on $C\times \Q$:
$$\O_{C\times \Q}^k \lra\widetilde
E\lra\widetilde E(D)\;,\leqno(4.25)$$
which, in turn,
induces an embedding $\Q\data \hookrightarrow\Q\tdata$, since
$\widetilde E(D)$ is a bundle of rank $r$
and degree $\tilde d$. What is less
obvious is the following.

\proclaim Theorem 4.26.
(i) For each $\delta >> 0$, there is a choice of
$\tau $ so that there is an embedding
$$\Q\data \hookrightarrow \B_\tau\tdata$$
where $\B_\tau\tdata$
is smooth of dimension $\tilde d k-r(k-r)(g-1)$ (see Theorem 3.20).
\hfil\break
(ii) If $\B = \B_\tau\tdata$
is chosen as in (i), and $U_\tau = U_\tau\tdata$ is the
universal rank $r$ bundle on
$C\times \B_\tau$, then the embedding in
(i) may be chosen so that the embedded $\Q\data$
is the scheme-theoretic intersection of $k\delta $
subvarieties, each of which is the zero-scheme of a map
$\O_\B \to U_\tau|_{\{p\}\times \B_\tau}$.

\noindent {\it Proof.} In Theorem 3.20, we showed that for
generic $\tau$ in the range given by Assumption 3.12, the moduli
space $\B_\tau\tdata$ is smooth, of the expected
dimension. Since $\B_\tau\tdata$ is a moduli space
for $\tau$-stable $k$-pairs, we need to show that the
family of $k$-pairs defined by the map (4.25) above is
$\tau$-stable for some $\tau$ in the desired range. But
for each $x\in \Q\data$, the map
$\O_C^k\to E(D) = \widetilde E(D)|_{C\times \{x\}}$
factors through $\O_C^k(D)$ and does not factor through any
subbundle of $E(D)$, so
if $E(D)$ fits in an exact sequence
$$0\lra E' \lra E(D) \lra E'' \lra 0\;,$$
then $\mu(E'') \geq \delta$
and  the desired $\tau$-stability follows if $\delta >> 0$.

Suppose that $D$ is the sum of $\delta >> 0$ distinct points
$p_1,...,p_\delta$ on $C$, and let
$\B_\tau = \B_\tau\tdata$ as above.
Then for each $p_i \in D$ and each
summand $e_j: \O_{\B_\tau} \hookrightarrow \O_{\B_\tau}^k$,
the universal
$k$-pair $\O_{C\times \B_\tau}^k \rightarrow U_\tau$ on $C\times
\B_\tau$ induces a section:
$\Phi_{i,j}: \O_{\B_\tau} \rightarrow U_\tau|_{p_i\times \B_\tau}.$
Let $Z_{i,j}\subset \B_\tau$ denote the zero-scheme of $\Phi_{i,j}$,
and
let $$Z = \bigcap _{1\leq i \leq \delta,1\leq j\leq k}Z_{i,j}\;.$$
If we restrict the universal $k$-pair to $C\times Z$, then the pair
factors:
$$\O_{C\times Z}^k \lra U_\tau(-D)|_{C\times Z}
\lra U_\tau|_{C\times Z}\;,\leqno(4.27)$$
and $\tau$-stability implies that
for each $z\in Z$, the restriction
of (4.27) to
$$\O_{C\times \{z\}}^k \to U_\tau(-D)|_{C\times
\{z\}}$$
cannot span any subbundle. Thus (4.27) determines a morphism
$Z\rightarrow \Q\data$ which inverts the map from $
\Q\data$ to $\B_\tau(\tilde d,r,k)$, and the theorem is proved.

Although the Quot scheme and Uhlenbeck compactification
are not in general
smooth, the following theorem and its corollary show that there is
at least some reasonable structure to these spaces for large degrees
relative to $g$, $r$, and $k$ (for a sharper description in the case
of $g=1$, see [Br]).

\proclaim Theorem 4.28.  There is a function $f(g,r,k)$ such that
$\Q\data$ is irreducible and generically reduced of the expected
dimension $kd-r(k-r)(g-1)$ for all $d\geq f(g,r,k)$.

\noindent {\it Proof.}  By induction
on the rank $r$ for fixed $g$ and $k$.  If $r=1$, then by Corollary
4.23, the Quot scheme $\Q\data$ is a smooth projective bundle of
dimension $kd-(k-1)(g-1)$ over $J_d$ as soon as $d\geq 2g-1$.  So
set $f(g,1,k)=2g-1$.  For rank two or greater we define
$$\eqalign{
\QS\data &=\{ (E,\vec\phi) : E \hbox{ is semistable and }\vec\phi
: \O_C^k\to E\cr
&\qquad\qquad\hbox{ is generically surjective }\}\;.\cr}$$
Then $\QS\data$ is  {\it simultaneously} an open subscheme of $
\Q\data$ and of $\B_\tau\data$ for any admissible $\tau$.  Since all
the $\B_\tau$'s are smooth and irreducible, of the expected
dimension, the theorem will follow once we show that every component
of the complement $\Q\data\setminus\QS\data$ has dimension smaller
than $kd-r(k-r)(g-1)$.

In order to bound the dimensions of this complement we first recall
that if $E$ is unstable then it fits into an exact sequence
$$0\lra S\lra E\lra Q\lra 0$$
where $S$ is stable.
Moreover, if $\vec\phi : \O^k_C\to E$ generically generates $E$ then
the induced section $\O_C^k\to Q$ must also generically generate.
Finally, if $d_s$, $r_s$ and $d_q$, $r_q$ denote the degrees and
ranks of $S$ and $Q$, respectively, because this exact sequence
exhibits the instability of $E$ we have $d_sr_q-d_qr_s > 0$.  Thus,
any $k$-pair $(E,\vec\phi)$ in the complement $\Q\setminus\QS$ may
be constructed as follows:
\itemitem{(a)} Choose non-negative integers $d_s, r_s$ and $d_q,
r_q$ such that $d_s+d_q=d$, $r_s, r_q > 0$, $r_s+r_q=r$, and
$d_sr_q-d_qr_s >0$.
\itemitem{(b)} Choose a stable bundle $S\in\N(d_s,r_s)$ and
$(Q,\vec\varphi)\in\Q(d_q, r_q, k)$.
\itemitem{(c)} Choose an extension $x\in\hbox{Ext}^1(Q,S)$ together
with a lift of the sections $\vec\varphi$ to sections $\vec\phi$ of
the bundle $E$ in the resulting exact sequence.

The main point is that the choice in (c) comes from a projective
space.  To be precise, if we let $V$ be defined by the long exact
sequence:
$$0\to\hbox{Hom}(Q,S)\to\hbox{Hom}(\O_C^k, S)\to V
\to \hbox{Ext}^1(Q,S)\to\hbox{Ext}^1(\O_C^k, S)\;,$$
then the choice in (c) is a point in $\P(V)$.  Furthermore, since
$S$ is stable with slope
at least $d/r$, we may assume that $\hbox{Ext}^1(\O_C, S)=0$, so we
have
$$\eqalign{\dim V &= k\chi(S)-\chi(Q^\ast\otimes S)\cr
&= k(d_s-r_s(g-1))-(d_s r_q - d_q r_s)+r_s r_q (g-1)\;.\cr}$$
The choices in (a) are discrete, so it suffices to count dimensions
for each choice of $d_s, r_s, d_q, r_q$.  We distinguish two cases.

\noindent {\it Case 1.}  Suppose $d_q \geq f(g, r_q, k)$.  Then
$\Q(d_q, r_q, k)$ has dimension $kd_q-r_q(k-r_q)(g-1)$ by induction,
so the dimension of the component coming from the choice of $d_s,
r_s, d_q, r_q$ is
$$\eqalign{\dim \Q(d_q, r_q, k) &+ \dim \N(d_s, r_s) +\dim \P(V)\cr
&= kd_q -r_q(k-r_q)(g-1) + r_s^2(g-1) +1 \cr
&\qquad k(d_s-r_s(g-1))-(d_s r_q- d_q r_s) + r_s r_q(g-1)-1\cr
&= kd -r(k-r)(g-1)-(d_s r_q-d_q r_s)- r_s r_q(g-1)\;,\cr}$$
which gives the desired bound, since we are assuming that $g\geq 1$.

\noindent {\it Case 2.}  Suppose $d_q < f(g, r_q, k)$.  Let
$f=f(g,r_q, k)$.  In this case, we do not know the dimension of
$\Q(d_q, r_q, k)$.  However, by the reasoning preceding Theorem
4.24, the dimension is bounded above by $kf-r_q(k-r_q)(g-1)$.
Hence, by the same calculation as in Case 1, we see that
$$\eqalign{\dim \Q(d_q, r_q, k) &+ \dim \N(d_s, r_s) +\dim \P(V)\cr
&\leq k(d_s+f)-r(k-r)(g-1) -(d_s r_q- d_q r_s) - r_s r_q(g-1)\cr
&= kd -r(k-r)(g-1)+k(f-r_q)-(d_s r_q-d_q r_s)- r_s r_q(g-1)\;.\cr}$$
But $k(f-r_q)\leq kf$ is bounded, and since we assumed that $d_q <
f$, the term $d_sr_q-d_qr_s$ grows with $d$.  For sufficiently large
$d$ we therefore have the desired inequality.  This completes the
proof.

\proclaim Corollary 4.29.  For $d \geq f(g,r,k)$ the Uhlenbeck
compactification $\U\data$ is irreducible and generically reduced.

\noindent {\it Proof.}
 From Theorem 4.24, the Uhlenbeck compactification is
the image of the Quot scheme, which is irreducible by Theorem 4.28.
Moreover, the two spaces share a dense open set, namely, $\M\data$,
so the Uhlenbeck compactification is generically reduced as well.

\beginsection 5.  Intersection numbers

\bigskip
\centerline{\it \S 5.1 Definitions}

\noindent
In this section, we rigorously define the intersection pairings
(1.2)  of the Introduction.  We show, in particular, that when the
degree $d$ is sufficiently large these pairings correspond to the
``definition" in (1.2) when $\M=\M\data$
is compactified by the Grothendieck Quot scheme of \S 4.2.  We will
use the following notation for intersections. If $c_1,\ldots, c_n$
are Chern classes of codimension $d_i$ on an irreducible projective
scheme $X$ such that $\sum_{i=1}^n d_i=\dim X$, then we will denote
by $\langle c_1\cdots c_n ; X\rangle$ the intersection pairing of
the $c_i$'s with $X$.  This is a well-defined integer, even if $X$
is not smooth.  We refer to Fulton [Ful] for details.

Recall that the evaluation map $\mu : C\times\M\to G(r,k)$ of
(1.1) defines Chern classes $X_1,\ldots , X_r$ on $\M$ by
pulling back the Chern classes of the tautological bundle $S^\ast$
on $G(r,k)$ and restricting to $\{p\}\times\M$.
Since the Quot scheme $\Q=\Q\data$ admits a universal rank $r$
bundle $\tilde E$ on $C\times\Q$ which extends the pullback of
$S^\ast$, we immediately obtain the following

\proclaim Lemma 5.1.  The Chern classes $c_i := c_i\left(\tilde E
\bigr|_{\{p\}\times \Q}\right)$ on $\Q$ restrict to the classes $X_1,
\ldots , X_r$ on $\M$.

Thus, as a first approximation one might expect the correct
definition of the intersection numbers (1.2) to be the following.
For any set of integers $s_1,\ldots, s_r$ such that $\sum_{i=1}^r
i s_i =\dim\Q$ one defines
$$\langle X_1^{s_1}\cdots X_r^{s_r}\rangle := \langle
c_1^{s_1}\cdots  c_r^{s_r} ; \Q\rangle \;.\leqno (5.2)$$
Unfortunately, there are problems with this definition.  It may be
the case that $\M$ has many components of different dimension, or
one component of dimension different than the expected dimension
$kd-r(k-r)(g-1)$ found in \S 3.  Or, the space $\M$ may have the
expected dimension but the ``compactification" given by the Quot
scheme (or the Uhlenbeck space) may contribute extra components,
perhaps even of the wrong dimension, which should not be counted in
the intersection numbers.

On the other hand, we showed in Theorem 4.28 that for sufficiently
large values of $d$ the Quot scheme $\Q\data$ is irreducible and
generically reduced of the expected dimension.  Thus, it follows that
(5.2) is a reasonable definition of the intersection numbers for
large $d$.
We can then use this fact to construct a good definition
for all $d$ as follows:

\proclaim Lemma 5.3.  If $D\subset C$ is an effective divisor of
degree $\delta$ chosen so that $\Q(d+r\delta, r,k)$ is irreducible
and generically reduced, let $\Q\data\hookrightarrow \Q(d+r\delta,
r,k)$ be the embedding defined in Theorem 4.26.  Then the pairing
$$\langle c_1^{s_1}\cdots c_{r-1}^{s_{r-1}}
c_r^{s_r + k\delta}; \Q(d+r\delta, r,k)\rangle\leqno(5.4)$$
is independent of  the choice of $D$.

\noindent {\it Proof.} Suppose $\Q(d)=\Q\data$ is irreducible of the
correct dimension, and $D=\{q\}$.  Then the restriction of the
universal bundle $\tilde E$ on $C\times \Q(d+r)$ to $C\times \Q(d)$
coincides with the twist $\tilde E(q)$.  Since the Chern classes
$c_1,\ldots, c_r$ are defined by restricting to $\{p\}\times \Q$,
the $c_i$'s extend without change from $\Q(d)$ to $\Q(d+r)$.  In
addition, using the universal sections $\vec\phi :\O^k\to \tilde E$,
we see as in Theorem 4.26 (i) that the image of $\Q(d)$ in $\Q(d+r)$
may be described as the intersection of the zero loci of the
sections $\phi_j : \O\to \tilde E_{\{q\}\times\Q(d+r)}$.  The zero
loci are {\it necessarily } regular because of the dimension of
$\Q(d)$.  This proves the lemma in this special case.  More
generally, suppose $\Q(d)$, $D$, and $D'$ are given, with
$\Q(d+r\delta)$ and $\Q(d+r\delta')$ irreducible and generically
reduced.  Then we may embed $\Q(d)$ in $\Q(d+r(\delta+\delta'))$ by
passing through either $\Q(d+r\delta)$ or $\Q(d+r\delta')$.  But the
special case of the lemma implies  that the intersection numbers in
both cases coincide with
$$\langle c_1^{s_1}\cdots c_{r-1}^{s_{r-1}}
c_r^{s_r+k(\delta+\delta')}; \Q(d+r(\delta+\delta')\rangle\;.$$
This completes the proof.

We therefore take (5.4) to be the definition of the intersection
numbers.  Next we will show that the intersection numbers may be
computed on the smooth moduli spaces $\B_\tau$ for certain choices
of $\tau$.  This will be essential for the computations in \S 5.3.

\proclaim Theorem 5.5.  Let $\Q(d)\hookrightarrow
\B_\tau(d+r\delta)$ be an embedding from Theorem 4.26 (i).  Let $c_i
:= c_i(U_\tau \bigr|_{\{p\}\times \B_\tau})$, where $U_\tau$ is the
universal rank $r$ bundle on $C\times \B_\tau$.  Then
$$\langle X_1^{s_1}\cdots X_r^{s_r}\rangle = \langle
c_1^{s_1}\cdots c_{r-1}^{s_{r-1}} c_r^{s_r+k\delta};
\B_\tau(d+r\delta)\rangle\;,\leqno(5.6)$$
where $\delta $ is as in Theorem 4.26.

\noindent {\it Proof.}  If $\Q(d)$ is irreducible of the expected
dimension, then as in the proof of Lemma 5.3 the $c_i$'s extend
unchanged to $\B_\tau$, each $Z_{i,j}$ in the intersection
$\Q(d)=\bigcap Z_{i,j}$ of Theorem 4.26 (ii) has codimension exactly
$r$,  and the intersection defines a subvariety of codimension
$rk\delta$ in $\B_\tau(d+r\delta)$.  Therefore, as in the proof of
Lemma 5.3, the pairing $\langle c_1^{s_1},\ldots, c_r^{s_r+k\delta};
\B_\tau(d+r\delta)\rangle$ is the same as the pairing
$\langle c_1^{s_1},\ldots, c_r^{s_r}; \Q(d)\rangle$.  The general
case now follows:  If $\Q(d)$ is fixed, let divisors $D$ and $D'$ be
chosen of degrees $\delta$ and $\delta'$, respectively, so that $D$
satisfies the conditions of Theorem 4.26 (i) and so that
$\Q(d+r\delta')$ is irreducible, of the correct dimension.  It then
follows from the definition of $\tau$-stability that both
$\B_\tau(d+r\delta)$ and $\Q(d+r\delta')$ embed in the same
space $\B_{\tau'}(d+r(\delta+\delta'))$, where
$\tau'=\tau+\delta'/r$.  But now as in the proof of Lemma 5.3, the
Quot scheme $\Q(d)$ embeds in
 $\B_{\tau'}(d+r(\delta+\delta'))$ either through $\B_\tau$ or
 through the Quot scheme, and by the special case of the lemma
 already proved, we see that both intersection numbers are computed
 by the same pairing of the $c_i$'s on
 $\B_{\tau'}(d+r(\delta+\delta'))$.

\bigskip
\centerline{\it \S 5.2 The conjecture of Vafa and Intriligator}

\noindent We now pause to describe a conjecture for the intersection
numbers defined above which is due to C. Vafa and was
worked out in detail
by K. Intriligator.  The conjecture arises from considerations of
certain superconformal field theories.  The physical reasoning which
led to this prediction may be found in [V], [I], and the references
therein;  here, however, we shall simply give the mathematical
formulation of the statement.

We begin by recalling that the ring structure of
$H^\ast\left(G(r,k), \C\right)$ is given by the free ring on the
Chern classes of the universal rank $r$ bundle $S$, modulo the ideal
of relations obtained by the vanishing of the Segre classes of $S$
beyond the rank of the universal quotient bundle $Q$.  More
precisely, let $X_i=c_i(S^\ast)$, $i=1,\ldots, r$, be as before.
 From the exact sequence (4.18) we have
$c_t(S^\ast)c_t(Q^\ast)=1$, where $c_t$ denotes the Chern
polynomial.  The fact that $Q^\ast$ has rank $k-r$ implies relations
on the $X_i$; if $I$
denotes the ideal generated by these relations,
then we have

\proclaim Proposition 5.7.  {\tenrm (cf. [B-T], p. 293)}  There is a
ring isomorphism
$$H^\ast\left(G(r,k),\C\right)\simeq \C[X_1,\ldots,X_r]/
 I\;.$$

Perhaps less well-known is the fact that $I$
is of the form $\langle
\partial W/\partial X_i ; i=1,\ldots, r
\rangle$, i.e. all the relations are obtained by setting to
zero the gradient of some homogeneous polynomial $W$ in the $X_i$'s.
To see this, write
$$c_t(Q^\ast)=\sum_{i=1}^k Y_i(X_1,\ldots,X_r) t^i\;.$$
Then $I$ is generated by the equations
$$Y_i=0\;,\; i=k-r+1,\ldots, k\;.\leqno(5.8)$$
On the other hand, making a formal expansion
$$-\log c_t(S^\ast)=\sum_{j\geq 0}W_j(X_1,\ldots, X_r)
t^j\;,$$
it is easily seen that
$$Y_{k+1-i}={\partial W_{k+1}\over \partial X_i}\;,\;
i=1,\ldots,r,$$
so that (5.8) implies the relations are generated by the equations
$dW_{k+1}=0$.
In terms of the Chern roots of $S^\ast$ defined by
$$c_t(S^\ast)=\prod_{i=1}^r(1+q_it)\;,$$
we find that we may take
$$ W=(-1)^{k+1} W_{k+1}=\sum_{i=1}^r{q_i^{k+1}\over k+1}\;.$$
We thus have

\proclaim Proposition 5.9.  Let $W$ be defined as above.  Then the
ideal $I$ in Proposition 5.7 is generated by the  polynomials
$\{\partial W/\partial X_i ; i=1,\ldots, r\}$.

In order to state the conjecture we first establish some notation.
 Let
$ W_1=W+(-1)^r X_1 $, and set
$$ h\xtuple
=(-1)^{r(r-1)/2}\det\left({\partial^2 W\over \partial X_i\partial
X_j }\right)\;.$$

\proclaim Conjecture 5.10. {\tenrm ([I], eq. (5.5))}
The intersection numbers (5.2) are given by
$$\langle X_1^{s_1}\cdots X_r^{s_r}\rangle
= \sum_{dW_1\ztuple=0} h^{g-1}\ztuple  Z_1^{s_1}\cdots Z_r^{s_r}\;.$$

\noindent We shall always interpret Conjecture 5.10 to apply to the
case where the degree $d$ is sufficiently large.

Consider, for example, the case $r=1$, i.e. the case of the maps to
projective space.  The polynomial $W=X^{k+1}/k+1$, and the critical
points of $W_1$ are the $k$-th roots of unity.  Furthermore,
$$h^{g-1}(Z)=\left(W^{\prime\prime}\right)^{g-1}=k^{g-1}
Z^{(k-1)(g-1)}\;.$$
Thus we have
$$\sum_{W_1^\prime(Z)=0} h^{g-1}(Z)Z^{kd-(k-1)(g-1)}=
\sum_{Z^k=1} k^{g-1} Z^{kd}
=\sum_{Z^k=1} k^{g-1}
= k^g\;.$$
By Theorem 2.9, $k^g$ is precisely the top intersection of the class
in the space of holomorphic maps of degree $d$ to $\P^{k-1}$ defined
by $X$.  We have established

\proclaim Theorem 5.11.  Conjecture 5.10 is true for $r=1$.

For the rest of the paper we shall assume $r=2$.  The reason for
this is that for maps into $G(2,k)$ the results of \S 3.3 and \S
5.1 give us an effective method for computing the left
hand side in Conjecture 5.10.
In the next subsection we shall set up this calculation for
arbitrary genus, however we shall only carry it through in the
case of elliptic curves $g=1$.
Therefore, let us set
$$\I =\sum_{dW_1(Z_1, Z_2)=0} Z_1^{kd-2n} Z_2^n\leqno(5.12)$$
for $n=0,1,\ldots, [kd/2]$
and $W_1$ the polynomial associated to
$G(2,k)$ as above.  Our goal is to simplify the expression (5.12).

As discussed in [I], the sum on the right hand side of (5.12)
is most clearly expressed in
terms of the Chern roots $q_1$ and $q_2$.  The critical points
of $W_1$  are given by
$$q_i =\alpha \xi_i\;,\;\alpha^k=-1\;,$$
where $\xi_1$ and $\xi_2$ run over the $k$-th roots of unity such that
$\xi_1\neq\xi_2$.  This overcounts
by a factor of 2, since the $X_i$'s are symmetric in the $q_i$'s.
Thus
$$\eqalign{
\I&= {(-1)^d\over 2}\sum_{\xi^k_i=1 \atop
\xi_1\neq\xi_2}(\xi_1+\xi_2)^{kd-2n}(\xi_1\xi_2)^n\cr
&={(-1)^d\over 2}\sum_{\xi_i^k=1}
(\xi_1+\xi_2)^{kd-2n}(\xi_1\xi_2)^n
-{(-1)^d\over 2}\sum_{\xi^k=1}(2\xi)^{kd-2n}\xi^{2n}\cr
&=
(-1)^{d+1} k 2^{kd-2n-1}
-{(-1)^{d+1}\over 2}\sum_{\xi_i^k=1}
(\xi_1+\xi_2)^{kd-2n}(\xi_1\xi_2)^n
\;.\cr}$$
If we now set $z=\xi_2\xi_1^{-1}$ we may eliminate one of the roots
and introduce a factor of $k$.  The result is
$$\leqalignno{ \I &=
(-1)^{d+1} k 2^{kd-2n-1}
-(-1)^{d+1}{k\over 2}\sum_{z^k=1} (1+z)^{kd-2n}z^n
\cr
&=(-1)^{d+1}k 2^{kd-2n-1}-
(-1)^{d+1}{k\over 2}\sum_{q=0}^{kd-2n}\sum_{z^k=1}
{kd-2n\choose q} z^{n+q}
 \cr
\I &=(-1)^{d+1}k 2^{kd-2n-1}
 - (-1)^{d+1}{k^2\over 2}\sum_{p\in{\Bbb Z}\atop n/k\leq p\leq d-n/k}
{kd-2n\choose kp-n}
\;,&(5.13)\cr}$$
since  the sum of $z^{n+q}$ over the $k$-th roots of unity vanishes
unless $n+q$ is of the form $kp$.  Equation (5.13) is the
desired expression.

\bigskip
\centerline{\it \S 5.3  Computations}

In this section we will outline a procedure for calculating all
intersection pairings of the form $\langle c_1^mc_2^n;
\B_\tau(d,2,k)\rangle$ where $\B_\tau$ is one of
the smooth moduli spaces of  $\tau$-stable $k$-pairs
from \S 3, and $m + 2n = kd - 2(k-2)(g-1) = \dim \B_\tau$.
 From these pairings and Theorem 5.5, one recovers all the
Gromov invariants for maps from a Riemann surface to $G(2,k)$. We
shall compute these invariants in
the case where $C$ is elliptic and show that they
agree with Conjecture  5.10.

For each fixed degree $d$,
recall that the admissible values of $\tau$ lie in the range
$d/2 < \tau < d-(2g-2)$ (Assumption 3.12)
and that if
$[\tau,\tau'] \cap {\Bbb Z} = \emptyset$ then $\B_\tau =
\B_{\tau'}$ (see the discussion following Theorem 3.38).
Thus, if we fix some $0 < \varepsilon < 1/2$, and list
moduli spaces: $\B_{d/2 + \varepsilon}, \B_{l \pm
\varepsilon}, [d/2] < l < d-(2g-2)$, then each smooth
moduli space appears in the list, and with the exception of
$\B_{d-(2g-1)+\varepsilon}$, each moduli space is counted
exactly twice; once as  $\B_{l-\varepsilon}$ and once as
$\B_{(l-1)+\varepsilon}$ (or $\B_{d/2+\varepsilon}$).

The calculation of the intersection pairings $\langle c_1^mc_2^n;
\B_\tau\rangle$ naturally breaks into two steps. First, an initial
calculation is necessary to compute the pairings
$\langle c_1^mc_2^n;\B_{d/2+\varepsilon}\rangle$.
If the degree is odd, then
by Theorem 3.26, the moduli space $\B_{d/2 + \varepsilon}$ is a
projective bundle over ${\cal N}(2,d)$, and the calculation is
analogous to the rank one calculation of Theorem 2.9.
In principle, this calculation is therefore always computable
by Grothendieck-Riemann-Roch
since the full cohomology ring
structure of ${\cal N}(2,d)$ is now known
(see [K], [T1], and [Z]). If the degree
is even, a further argument is necessary because of the existence of
semistable bundles which are not stable. This can be dealt
with  without too much difficulty in the genus one case, as we
will see below.

Next, there are the flip calculations which compute
the {\it difference}
between the pairings $\langle c_1^mc_2^n;\B_{l-\varepsilon}\rangle$
and $\langle c_1^mc_2^n;\B_{l+\varepsilon}\rangle$. This
is represented by a class which lives on a ``rank one''
locus $\pspace(d-l,k) \times  J_l$,
and thus may be viewed as a rank one calculation.
Since every smooth moduli space is a $\B_{l \pm \varepsilon}$,
we see that all the pairings, including those which produce
the Gromov invariants, are
computed by the initial calculation followed by a sequence of flip
calculations. We now give the details of these two
contributions.

\bigskip
\noindent
{\it The Initial Calculations:} Because the moduli spaces $
\B_{d/2+\varepsilon}$ behave differently depending upon the parity of
$d$, we will consider the even and odd cases separately.

\bigskip
\noindent {\it Odd degree:} If $V$
is any universal bundle on $C\times {\cal N}(2,d)$, then
by Theorem 3.26 the moduli space $\B_{d/2+\varepsilon}$ is
isomorphic to the projective bundle $\P\left((\rho_*V)^{\oplus
k}\right)$,
and if $\pi:\B_{d/2+\varepsilon} \rightarrow {\cal N}(2,d)$ is the
projection, then by Corollary 3.32 the universal bundle on $C\times
\B_{d/2+\varepsilon}$ is isomorphic to $\pi^*V \otimes {\cal O}(1)$ (all
pulled back to $C\times \B_{d/2+\varepsilon}$). The intersection
pairings now reduce
to a calculation involving Grothendieck-Riemann-Roch.
In particular, we easily  obtain

\proclaim  Calculation 5.14.
If $C$ is elliptic and $d$ is odd, then
$$\langle c_1^mc_2^n;\B_{d/2+\varepsilon} \rangle = k2^{m-1}.$$

\noindent
{\it Proof.} The moduli space ${\cal N}(2,d)$ is isomorphic to $C$
itself in this case.  Moreover, if we let $\Delta \subset
C\times C$ be the diagonal, then
there is a canonical extension on $C\times {\cal N}(2,1)$:
$$0 \rightarrow \O_{C\times C} \rightarrow V \rightarrow
\O_{C\times C}(\Delta) \rightarrow 0$$
which produces a universal bundle. If we let $D\subset C$ be a
divisor of degree $(d-1)/2$, then $V(D):= V\otimes \rho^*
O_C(D)$ is a universal bundle
on $C\times {\cal N}(2,d)$, which satisfies:
$$\eqalign{ (i)\quad  c_1 (V(D)_p) &= [q] \cr
(ii) \quad  c_1 (\rho_* V(D)) &= ((d-1)/2) [q]\;,\cr} $$
where $[q]$ is the class of a point in ${\cal N}(2,d)$.

We need to consider the bundle $(\pi^*V(D) \otimes \O(1))_p$. If
we let $z$ represent the first Chern class of $\O (1)$, then
from (i) and the formula for the Chern classes of a tensor product,
$c_1 =  \pi^*[q] + 2z$ and $c_2 = \pi^*[q]z + z^2 = c_1^2/4$.
Thus $c_1^mc_2^n =  2^{-2n} c_1^{kd}$. If we let
$s_0$ and $s_1$ be the first two Segre classes of
$(\rho_*V(D))^{\oplus k}$,
then the intersections are calculated as follows (cf. (2.12)):
$$\eqalign{
c_1^{kd}&= (\pi^*[q] + 2z)^{kd} =
(kd)2^{kd-1}(s_0)\pi^*[q] + 2^{kd}(s_1)\cr
&= (kd)2^{kd-1} - 2^{kd}k\left({d-1\over 2}\right) = k2^{kd-1}
\;,\cr}$$
and the lemma is proved.

\bigskip
\noindent {\it Even degree:} Here we unfortunately do not have a
simple description of the moduli spaces $\B_{d/2 +
\varepsilon}$ in general. In case the genus is one, however,
the following proposition will enable us to compute the intersection
numbers on a projective variety mapping with degree two
onto the moduli  space.

Suppose $C$ is elliptic, and
let $P=P(d/2,k) =
\pspace(d/2,k)$ be the rank one moduli space
defined in (2.4).
Let $\Delta_P \subset  P\times P$ be
the diagonal, and let $ P(2)$ be the blow-up of $ P\times
 P$ along the diagonal. Let $A\subset  P(2)$ be the
exceptional divisor, and let ${\cal L}(2)$ be the pull-back
to $C\times  P(2)$ of the
direct sum of the two obvious universal line bundles on $C\times
 P(2)$. Finally, let ${\cal L}_A$ be the pull-back
to $C\times A$ of the universal bundle on
$C\times \Delta _P \cong C \times  P$. Then:

\proclaim  Lemma 5.15. (i) There is a holomorphic map
$\sigma: P(2) \rightarrow \B_{d/2+\varepsilon}$,
whose fibres are generically finite of degree two.
(ii) There is an exact sequence of sheaves on $C\times  P(2)$:
$$0 \rightarrow \sigma^*U_{d/2+\varepsilon} \rightarrow {\cal L}(2)
\rightarrow {\cal L}_A \rightarrow 0$$

\noindent {\it Proof.}
The diagonal inclusion of line bundles
${\cal L} \hookrightarrow {\cal L}\oplus {\cal L}$
over $C\times\Delta_P$ induces the exact sequence of bundles on
$C\times A$:
$$0 \rightarrow {\cal L}_A \rightarrow {\cal L}(2)|_A \rightarrow
{\cal L}_A \rightarrow 0$$
Denote by ${\cal K}$ the kernel of the
sheaf map ${\cal L}(2) \rightarrow {\cal L}_A$, where ${\cal L}_A$
is extended by zero to all of $C\times  P(2)$. Then the $k$
sections ${\cal O}^k \rightarrow {\cal L}(2)$ lift to sections of
${\cal K}$, so we get in this way a family of pairs parametrized by
$ P(2)$. One checks that these pairs are all $d/2 + \epsilon$-
stable, which gives the map $\sigma $ of (i). But this also gives
(ii), since $\B_{d/2+\epsilon}$ is a fine moduli space.

The following proposition will be essential to both this calculation
and the flip calculations which follow.

\proclaim  Proposition 5.16. Suppose $X\subset Y$ is a
proper inclusion of smooth,
projective varieties of dimensions $M$ and $N$, respectively.
Suppose that $E$ is a vector bundle of rank two on $Y$, such that
when $E$ is restricted to $X$, there is a surjective map
$E|_X \rightarrow L$ to an invertible sheaf $L$ on $X$.
Let $\tilde Y$ denote the blow-up of $Y$ along $X$, let
$A\subset \tilde Y$ denote the exceptional divisor, and let $\tilde E$
and $\tilde L$ denote the pullbacks of $E$ and $L$ to $\tilde Y$ and
$A$, respectively. Let $i:A\hookrightarrow \tilde Y$ be the
inclusion map, and finally, let $F$ be the rank two bundle defined as
the kernel in the sequence:
$$0 \rightarrow F \rightarrow \tilde E \rightarrow i_*\tilde L
\rightarrow 0.$$
Then:
\itemitem{(i)} $c_1(F) = c_1(\tilde E) - A$ and $c_2(F) =
c_2(\tilde E) - i_*(c_1(\tilde L))$.
\itemitem{(ii)} If $m,n$ are integers such that $m+2n = N$, then
the integer
$(c_1^m(F)c_2^n(F) - c_1^m(\tilde E)c_2^n(\tilde E);\tilde Y)$ is
equal to the evaluation on $X$ of the coefficient of $t^{M-n}$ in
the polynomial:
$$-s_t(\nu(X))(1+c_1(E)t)^m(c_1(L)+c_2(E)t)^n$$
where $s_t(\nu(X))$ is the Segre polynomial of the normal bundle to
$X$ in $Y$.

\noindent {\it Proof.} (i) is an immediate consequence of the
Grothendieck-Riemann-Roch theorem applied to the sequence defining
$F$ (see [Ful]).
To prove (ii), note first that because of (i), we have
$$c_1^m(F)c_2^n(F)  =
\left(\sum _{a=0}^m \left({m \atop a}\right)
(-A)^ac_1^{m-a}(\tilde E)\right)
\left(\sum _{b=0}^n \left({n \atop b}\right)
(-1)^b(i_\ast c_1({\tilde L}))^bc_2^{n-b}(\tilde E)\right).$$
When the right hand side is evaluated on $\tilde Y$, all the
terms with $a=0$ but $b\ne 0$ vanish. Since all the
other terms may be evaluated on $A$, the right hand side becomes:
$$c_1^m(\tilde E)c_2^n(\tilde E)  +
\sum _{a=1}^m \sum _{b=0}^n
\left({m \atop a}\right) \left({n \atop b}\right)
(-1)(-A)^{a+b-1}c_1^{m-a}(\tilde E)c_2^{n-b}(\tilde E)
c_1^b({\tilde L})$$
where $c_1^m(\tilde E)c_2^n(\tilde E)$ is evaluated on
$\tilde Y$ and the double sum is evaluated on $A$.
Since $m+2n = N$, see that the
restriction of $(-A)^{a+b-1}$ to $A$
is realized by the Segre class
$s_{a+b+M-N}(\nu (X))$
whenever $a+b \ge N-M$. Since the restriction to $A$ of the other
classes $c_1(\tilde E),c_2(\tilde E)$ and $c_1(\tilde L)$
are all pullbacks from $X$, it follows that the
terms involving lower powers of $-A$ all evaluate to zero on $A$.
Thus we can substitute the appropriate Segre classes for powers of
$-A$ and part (ii) immediately follows.

\proclaim  Calculation 5.17. If $C$ is elliptic and $d$ is
even, then:
$$\langle c_1^mc_2^n;\B_{d/2+\varepsilon}\rangle =
{k^2\over 2} { kd-2n \choose kd/2 - n}
- k2^{m-1}.$$

\noindent {\it Proof.} Let $\sigma$ be the degree two map of Lemma 5.15,
and,  abusing notation slightly, let
$c_1 = \sigma^*(c_1)$ and $c_2 = \sigma^*(c_2)$.
Then part (i) of Lemma 5.15 implies that the evaluation of
$c_1^mc_2^n$ on $P(2)$ gives twice the evaluation on
$\B_{d/2+\varepsilon}$.

If we apply Proposition 5.16 to the exact sequence in part (ii) of
Lemma 5.15, then we may express $\langle c_1^mc_2^n; P(2)\rangle$
as the sum of two contributions, the first being
$$\langle c_1^m({\cal L}(2)_p) c_2^n({\cal L}(2)_p); P(2)\rangle$$
and the second being the
evaluation on $\Delta_P$ of the coefficient of
$t^{m/2}$
in the polynomial
$$-s_t(\nu(\Delta_P))(1+c_1({\cal
L}(2)_p)t)^m(c_1({\cal L}) + c_2({\cal L}(2))t)^n.$$
The first contribution is easy to compute. If we let ${\cal L}_1$
and ${\cal L}_2$ be the two pullbacks of ${\cal L}_p$ to $
P\times  P$, then $c_1({\cal L}(2)_p) = c_1({\cal L}_1) +
c_1({\cal L}_2)$ and $c_2({\cal L}(2)_p) = c_1({\cal L}_1)c_1({\cal
L}_2)$. Then:
$$\eqalign{
\langle c_1^m({\cal L}(2)_p)c_2^n({\cal L}(2)_p); P(2)\rangle
&=
\langle (c_1({\cal L}_1)+c_1({\cal L}_2))^mc_1({\cal L}_1)^nc_1({\cal
L}_2)^n; P\times P\rangle\cr
&= { m \choose m/2} \langle c_1^{kd/2}({\cal L}_1)
c_1^{kd/2}({\cal L}_2); P\times  P\rangle\cr
&= {kd-2n\choose  kd/2 - n}k^2\;,\cr}$$
since $c_1^{kd/2}({\cal L}_i) = k$ by Theorem 2.9. This gives
the first term in the calculation.
For the second contribution, observe that if $L_{d/2}$
is a universal line
bundle on $C\times  J_{d/2}$, and $T P$ is the
tangent bundle to $P$, then
$$s_t(\nu(\Delta_P)) = s_t(T P) =
s_t((\rho_*L_{d/2})^{\oplus k}(1))$$
from the relative Euler sequence and the fact that $TC =
\O_C$.

Since in addition, we have $c_1({\cal L}(2))|_{\Delta _P} =
2c_1({\cal L})$ and $c_2({\cal L}(2))|_{\Delta _P} =
c_1^2({\cal L})$, we see that the second contribution is one-half the
evaluation on $P$ of the $t^{m/2}$ coefficient of
$$- c_t(\oplus^k(\rho_*L_{d/2})(1))^{-1}
(1+2c_1({\cal L}_p)t)^m
(c_1({\cal L}_p) + c_1^2({\cal L}_p)t)^n\;.\leqno(5.18)$$
Let $a \in {\Bbb Z}$ be defined by $a[\Theta] = c_1(L_{d/2})_p$ and
let $[z] = c_1(\O_P(1))$. Then recall that $c_1({\cal
L}_p) = a[\Theta] + [z]$ and this is independent of the choice of
$L_{d/2}$. It follows from Grothendieck Riemann-Roch that
$c_1(\oplus ^k(\rho_*L_{d/2})) = k(ad/2 - 1)[\Theta]$.
The computations become more manageable if we make the following
change of variables:
$$w = z + (a- 2/d)\Theta, \ \ \varepsilon =  (2/d) \Theta\;.$$
This change of variables is motivated by the observation that
if there were a choice of universal bundle $L_{d/2}$ such that
$a = 2/d$, then $\rho_*{L_{d/2}}$ would have the
Chern classes of the trivial bundle, and $c_t((\rho_*
L_{d/2})^{\oplus k}(1))$
would be simply $(1+zt)^{kd/2}$. Such a universal
bundle does not exist for $d > 2$  since $2/d$ is not an integer;
nevertheless, via this change of variables we may pretend that $w$
represents $\O(1)$ for such a universal line bundle.
In terms of the new variables, we have the following identities:
$$c_t((\rho_*L_{d/2})^{\oplus k}(1)) = (1+wt)^{kd/2},
\ \ \varepsilon^2=0, \ \ w^{kd/2}=0, \ \ w^{kd/2 -1}\varepsilon = 2/d
\;.$$
Then we  have
$$\eqalign{
(5.18) \ \ &= -
{(1+2(w+\varepsilon)t)^m(w+\varepsilon)^n(1+(w+\varepsilon)t)^n
\over
(1+wt)^{kd/2}}\cr
&=  - w^n \left( {(1+2(w+\varepsilon)t)^m\over (1+wt)^{m/2}}\right)
- n\varepsilon w^{n-1} \left(  {(1+2wt)^{m+1}\over(1+wt)^{m/2+1}}
\right)\cr
&= - 2m\varepsilon w^nt
\left(  {(1+2wt)^{m-1}\over(1+wt)^{m/2}} \right)
- n\varepsilon w^{n-1} \left(  {(1+2wt)^{m+1}\over(1+wt)^{m/2+1}}
\right)\;.\cr}$$
 From the identity: $(1+2st)^{2r-1}(1+st)^{-r} =
2^{2r-2}s^{r-1} + $ other powers of $t$, we see that
the coefficient of $t^{m/2}$ in $(5.18)$ evaluates to
$- 2m( 2/d)2^{m-2} - n(2/d)2^m = - k2^m$ on $P$,
which finishes the calculation.

\bigskip
\noindent {\it The Flip Calculations:}
Recall the ``flip'' diagram of Theorem
3.42:
$$\matrix{&& \tilde\B_l &&\cr
&{}^{p_-}\swarrow &&\searrow {}^{p_+}&\cr
\B_{l-\varepsilon}&&&& \B_{l+\varepsilon}\cr}$$
and its restriction to the exceptional divisor
$$\matrix{&& A &&\cr
&\swarrow &&\searrow &\cr
\P(W^-)&&&&\P(W^+)\cr
&{}_{\sigma_-}\searrow &&\swarrow_{\sigma_+}&\cr
&& P(d-l,k)\times J_l &&\cr}$$
For the duration of the calculations it will be convenient to adopt
the following conventions for naming various bundles:
\itemitem{(1)} Denote
by $U^\pm$ the universal rank two bundles on $C\times
\B_{l\pm \varepsilon}$. Denote by $U^\pm_p$ the restriction of
$U^\pm$ to $\{ p \} \times \B_{l\pm \varepsilon}$, and finally
denote by $c_i^\pm$ the pull-back of the
$i$th ($i=1,2$) Chern class of $U^\pm_p$ to $\tilde \B_l$.

\noindent Thus,
according to this notation the flip contribution is the
evaluation on the fundamental class of $\tilde\B_l$ of the
polynomial $(c_1^-)^m(c_2^-)^n - (c_1^+)^m(c_2^+)^n$.

\itemitem{(2)} Denote by $L_d$ any choice of universal line bundle on
$C\times  J_d$.
Denote by $L_d(x)$ the universal line bundle
on the moduli space $ P(d,k)$,  so if $ P(d,k)$
is identified with the projective bundle
$\P( (\rho_* L_d)^{\oplus k})$,
then $\O(x) = \O_P(1)$.

\itemitem{(3)} Denote by $\O_{\P(W^-)}(z)$ the line bundle
$\O_{\P(W^-)}(1)$, so by Proposition 3.36, the normal
bundle $\nu(\P(W^-)) = W^+(-z)$.

Throughout the computations, if $E$ is a vector bundle on a
variety $Y$ and $f:X\rightarrow Y$ is understood, then we will
denote also by $E$ the pullback of $E$ to $Y$. We trust that the
possible confusion arising from this is less than the confusion that
would be caused by the proliferation of notation necessary to make
everything precise. For example, we have already used
$\O_P(1)$ to denote its pullback to $C\times
P$ and $W^+$ to denote its pullback to ${\P}(W^-)$.

\proclaim Proposition 5.19. There is an exact sequence of
sheaves on $C\times \tilde \B_l$:
$$ 0 \rightarrow U^-(-A) \rightarrow U^+ \rightarrow L_{d-l}(x)
\rightarrow 0$$
where $L_{d-l}(x)$ is extended by zero from $C\times A$.

\noindent {\it Proof.} Let ${\cal K}$ be the kernel of the
map $U^+ \rightarrow L_{d-l}(x)$ obtained by pulling back the map
in Proposition 3.44 (ii).
The $k$ sections $\O^k
\rightarrow U^+$ pull back after twisting to give
a map $\O^k(-A) \rightarrow {\cal K}$, hence $k$ sections
of ${\cal K}(A)$. If this results in a family of $\tau - \varepsilon$
stable pairs, then the result follows immediately from the universal
property of blowing up. But this is precisely the content of (3.14)
of Thaddeus [T2] in the case where $k=1$, and the
general case is the same.

The following global version of Proposition 3.39 will
be essential to calculate the Segre polynomial of
$\nu(\P(W^-))$ (see also [T2]):
If $\rho: C\times  P(d-l,k) \times  J_l \rightarrow
 P(d-l,k) \times  J_l$ is
the projection, then we may write:
\itemitem{(i)} $W^- = R^1\rho_*(L_l^*\otimes L_{d-l}(x))$ and
\itemitem{(ii)} $0 \rightarrow \rho_*(L_l\otimes  L_{d-l}^*(-x))
\rightarrow (\rho_\ast L_l)^{\oplus k} \rightarrow W^+ \rightarrow
R^1\rho_*(L_l\otimes L_{d-l}^*(-x))\rightarrow 0$.

 From now on, assume that $C$ is elliptic. Then the long exact
sequence in (ii) reduces to a short exact sequence, since
$R^1\rho_*(L_l\otimes L_{d-l}^*(-x)) = 0$ if $2l-d > 0$. This
also implies that the first term of (ii) is a vector bundle, of rank
$2l-d$, and we will need the following:

\proclaim Lemma 5.20. The top
Chern class $c_{2l-d}\left(\rho_*\left(L_l\otimes L_{d-l}^*(-x)\right)
\otimes \O_{\P(W^-)}(-z))\right)$ vanishes.

\noindent {\it Proof.} Consider the following general setup.
Suppose $X$ is any smooth projective variety,
${\cal L}$ is a line bundle on $C\times X$, and $\rho :
C\times X \rightarrow X$ is the projection. Let $E = \rho_*{\cal L}$
and let $F = R^1\rho_*{\cal L}^*$, and suppose that $R^1\rho_*{\cal L}$
and $\rho_*{\cal L}^*$ both vanish, so $E$ and $F$ are vector bundles
on $X$, of the same rank n, since $C$ is assumed to be elliptic.
Now, by Grothendieck-Riemann-Roch and the fact that the Todd class
of an elliptic curve is trivial, it follows that the Chern classes
of $E$ are the same as the Chern classes of $F^*$, so since the top
Chern class $c_n(F^*\otimes \O_{{\P}(F)}(-1))$ is easily
seen to vanish, it follows that $c_n(E\otimes \O_{\P
(F)}(-1))$ vanishes as well. The lemma follows if we let
${\cal L}$ be $L_l\otimes L_{d-l}^*(-x)$.

Armed with 5.19 and 5.20, we can finally compute:

\proclaim  Calculation 5.21. If $C$ is elliptic and
$l \in (d/2,d)$, then we have:
$$\langle c_1^mc_2^n;\B_{l-\varepsilon}(d,2,k))\rangle-
\langle c_1^mc_2^n;\B_{l + \varepsilon}(d,2,k)\rangle
= (-1)^d k^2 \left( {kd-2n} \atop {kl-n} \right)\;,$$
where by convention the right hand side is zero if $kl-n <0$ or
$kl-n > kd-2n$.

\noindent
{\it Proof.}
By Propositions 5.16 and 5.19, we see that the quantity we need to
compute is the evaluation on $\P(W^-)$
of the coefficient of $t^{\dim\P(W^-)-n}$ in the
power series:
$$ -c_t(\nu(W^+(-z))^{-1}(1+c_1^-t)^m
(c_1(L_{d-l}(x))_p + c_2^-t)^n\;.\leqno(5.22)$$
 From the exact sequence (ii) above, we have:
$$c_t(\nu(W^+(-z))^{-1} =
c_t(\oplus^k(\rho_*L_l)(-z))^{-1}
c_t(\rho_*(L_l\otimes L^*_{d-l}(-x))(-z))$$
and from part (i) of Proposition 3.44, we see that when restricted to
$\P(W^-)$, we have:
$$\eqalign{
c_1^- &= c_1(L_{d-l}(x))_p + c_1(L_l)_p - z \cr
c_2^- &= c_1(L_{d-l}(x))_p (c_1(L_l)_p -z)\;.\cr}$$
Suppose $c_1(L_l)_p = a[\Theta _l]$ and $c_1(L_{d-l})_p = b[\Theta
_{d-l}]$. Then as in Calculation 5.17, the following change
of variables simplifies things considerably:
$$\eqalign{
\varepsilon_l &= {1\over l}\Theta_l\cr
y &= x + b\Theta_{d-l}+\varepsilon_l \cr}
\qquad\qquad
\eqalign{
\varepsilon_{d-l} &= {1\over d-l}\Theta_{d-l}\cr
w &= z - a\Theta_l  +  \varepsilon_l\;,\cr}$$
and while we're at it, let $E = \rho_*(L_l\otimes L^*_{d-l}(-x))(-z)$.

In terms of these variables, the power series (5.22) becomes:
$$\eqalign{(5.22)&=
-(1-wt)^{-kl}c_t(E)(1-wt+yt)^m(y - \varepsilon_l)^n
(1-wt+\varepsilon_lt)^n\cr
&= -c_t(E)\left[ \sum_{p=0}^m { m \choose p }
(yt)^p(1-wt)^{N-p}
\left\{ y^n + n\varepsilon_l ty^n(1-wt)^{-1}
-n\varepsilon_l y^{n-1}\right\}\right]\;,\cr}$$
where $N=n+m-kl$, since $\varepsilon_l^2 = 0$.
In addition, one checks that when evaluated on $P
(d-l,k)\times  J_l$, the following hold:
$\varepsilon_ly^{k(d-l)} =  k/l$, $y^{k(d-l)+1} =
(k(d-l)+1)k/l$ and $y^p = 0$ for $p > k(d-l)+1$.
These, together with Lemma 5.20, imply that only three terms
contribute to the coefficient of $t^{{\rm dim}(\P(W^-)) - n}$;
namely:
$$\eqalign{
-c_t(E) t^{-n} \Biggl\{ \left( m \atop {k(d-l)-n+1}\right)
&\left[  {(yt)^{k(d-l)+1}\over 1-wt} -
n\varepsilon_l t {(yt)^{k(d-l)}\over 1-wt}\right] \cr
+&\left( m \atop {k(d-l)-n}\right) n\varepsilon _l t
{(yt)^{k(d-l)}\over 1-wt} \Biggr\}\;.\cr}$$
Now, since $c_t(E)$ behaves like $(1-wt)^{2l-d}$ when paired with
the classes $y^{k(d-l)+1}$ and $\varepsilon_ly^{k(d-l)}$,
the calculation is reduced to finding the coefficient of
$t^{2l-d+k(d-l)-n}$ in the following polynomial:
$$\eqalign{
-(1-wt)^{2l-d-1}t^{k(d-l)+1-n}
\biggl\{ &\left( m \atop {k(d-l)-n+1}\right)
\left[ y^{k(d-l)+1} - n\varepsilon_l y^{k(d-l)} \right] \cr
& + \left( m \atop {k(d-l)-n}\right) n\varepsilon _l
y^{k(d-l)}\biggr\} \;,\cr}$$
which evaluates to
$$-(-1)^{2l-d-1}
\left\{ \left( m \atop {k(d-l)-n+1}\right)
{ (k(d-l)-n+1)k\over l} +
\left( m \atop {k(d-l)-n}\right)
 {nk \over l} \right\}$$
$$ = (-1)^d k^2 \left( {kd-2n}\atop {kl-n}
\right)$$
as desired.

Finally, we use calculations 5.14, 5.17 and 5.21 to compute the
Gromov invariants associated to maps from an elliptic curve $C$ to
$G(2,k)$. Suppose $d$ is a nonnegative integer. Then the expected
dimension of $\Q(d,2,k)$ is $dk$, and thus by
Theorem 5.5, the intersection pairing  $\langle X_1^mX_2^n
\rangle$ for $m+2n=kd$
is realized as $\langle c_1^mc_2^{n+k\delta};\B_\tau (d+2\delta)
\rangle$ for
sufficiently large $\delta$. Also, $\tau $, following the proof
of Theorem 4.24, may be
chosen to be $d +\delta + \varepsilon$.

This pairing on $\B_\tau (d+2\delta)$ is now computed by first
computing the corresponding pairing on $\B_{d/2+\varepsilon}$,
which is either $k2^{m-1}$ if $d$ is odd, by Calculation 5.14, or
$-k2^{m-1} +  {k^2\over 2} {m \choose m/2} $ if $d$ is
even, by Calculation 5.17. The correction terms, which measure how
the pairing changes as $\tau $ moves from $d/2 + \delta +
\varepsilon$ to
$d + \delta + \varepsilon$, are in both the even and odd case equal to
$$(-1)^d k^2 \sum _{l= [d/2+1]+\delta}^{d+\delta} { m \choose
kl - (n+k\delta)}  =
(-1)^d k^2 \sum _{l= [d/2+1]}^{d} { m \choose
kl - n} $$
by repeated application of Calculation 5.21.
If we combine these results, we finally get:
$$\langle X_1^{kd-2n} X_2^n\rangle =(-1)^{d+1}k 2^{kd-2n-1}
 - (-1)^{d+1}{k^2\over 2}\sum_{p\in{\Bbb Z}\atop n/k\leq p\leq d-n/k}
{kd-2n\choose kp-n}\;,$$
which proves Theorem 1.7.
Thus, the conjecture of Vafa and Intriligator is true in these
cases as well:

\proclaim Theorem 5.27.  Conjecture 5.10 is true for maps from an
elliptic curve to $G(2,k)$.

\beginsection {} References

\item{[A-B]} Atiyah, M. F. and R. Bott, {\it The Yang-Mills equations
over Riemann surfaces}, Phil. Trans. R. Soc. Lond. A {\bf 308}
 (1982), 523-615.

\item{[A-C-G-H]}  Arbarello, E., M. Cornalba, P. Griffiths, and J.
Harris, ``Geometry of Algebraic Curves", Vol. I,  Springer-Verlag, New
York, Berlin, Heidelberg, 1985.

\item{[B]} Bradlow, S. B., {\it Special metrics and stability for
holomorphic bundles with global sections}, J. Diff. Geom. {\bf 33}
(1991), 169-214.

\item{[B-D]} Bradlow, S. B. and G. D. Daskalopoulos, {\it Moduli of
stable pairs for holomorphic bundles over Riemann surfaces}, Int.
J. Math. {\bf 2} (1991), 477-513.

\item{[B-D-W]} Bradlow, S., G. Daskalopoulos, and R. Wentworth, {\it
Birational equivalences of vortex moduli}, preprint, 1993.

\item{[B-T]} Bott, R. and L. Tu, ``Differential Forms in Algebraic
Topology", Springer-Verlag, New York, Berlin, Heidelberg, 1982.

\item{[Be]} Bertram, A., {\it Moduli of rank 2 vector bundles,
theta divisors, and the geometry of curves on projective space},  J.
Diff. Geom. {\bf 35} (1992), 429-469.

\item{[Br]} Brugui\`eres, A., {\it The scheme of morphisms from an
elliptic curve to a Grassmannian}, Comp. Math. {\bf 63} (1987),
15-40.

\item{[D]} Daskalopoulos, G.D., {\it The topology of the space of
stable bundles over a compact Riemann surface}, J. Diff. Geom. {\bf
36} (1992), 699-746.

\item{[F]} Floer, A., {\it Symplectic fixed points and holomorphic
spheres}, Commun. Math. Phys. {\bf 120} (1989), 575-611.

\item{[Ful]} Fulton, W., ``Intersection Theory",  Springer-Verlag,
New York, Berlin, Heidelberg, 1984.

\item{[G1]} Gromov, M.,  {\it Pseudo holomorphic curves in
symplectic manifolds}, Invent. Math. {\bf 82} (1985), 307-347.

\item{[G2]} Gromov, M.,  {\it Soft and hard symplectic geometry}, in
Proceedings of the International Congress of Mathematicians,
Berkeley, 1986.

\item{[G-H]} Griffiths, P. and J. Harris, ``Principles of Algebraic
Geometry", Wiley, New York, 1978.

\item{[G-S]} Guillemin, V. and  S. Sternberg, {\it Birational
equivalence in the symplectic category}, Invent. Math. {\bf 97}
(1989), 485-522.

\item{[Gro]} Grothendieck, A.,  {\it Techniques de construction et
the\'eor\`emes d'existence en g\'eometrie alg\'ebrique IV: Les
sch\'emas de Hilbert}, S\'eminaire Bourbaki {\bf 221} (1960/61).

\item{[I]} Intriligator, K., {\it Fusion residues}, preprint, 1991.

\item{[K2]} Kirwan, F., {\it On the homology of compactifications of
moduli spaces of vector bundles over a Riemann surface}, Proc.
London Math. Soc. (3) {\bf 53} (1986), 237-266.

\item{[N]} Newstead, P.E., ``Introduction to Moduli Problems and
Orbit Spaces", Tata Inst. Lectures {\bf 51}, Springer-Verlag,
Heidelberg, 1978.

\item{[R]} Ruan, Y., {\it Toplogical sigma model and Donaldson type
invariants in Gromov theory}, preprint, 1993.

\item{[S-U]} Sacks, J., and K. K. Uhlenbeck, {\it The existence of
minimal immersions of 2-spheres}, Ann. Math. {\bf 113} (1981),
1-24.

\item{[Str]} Str\o mme, S., {\it On parameterized rational curves in
Grassmann varieties},  in Lecture Notes in Math. {\bf 1266},
Springer-Verlag, New York, Berlin, Heidelberg, 1987.

\item{[T1]} Thaddeus, M., {\it Conformal field theory and the
cohomology of the moduli space of stable bundles}, J. Diff. Geom.
{\bf 35} (1992), 131-149.

\item{[T2]} Thaddeus, M., {\it Stable pairs, linear systems, and the
Verlinde formula}, preprint, 1992.

\item{[Ti]} Tiwari, S.,  preprint.

\item{[V]} Vafa, C.,  {\it Topological mirrors and quantum rings},
in ``Essays on Mirror Manifolds", S.-T. Yau, ed., International
Press, Hong Kong, 1992.

\item{[Wf]} Wolfson, J., {\it Gromov's compactness of
pseudo-holomorphic curves and symplectic geometry}, J. Diff. Geom.
{\bf 28} (1988), 383-405.

\item{[Wi]} Witten, E., {\it Topological sigma models}, Commun.
Math. Phys. {\bf 118} (1988), 411-449.

\item{[Z]} Zagier, D., unpublished.

\bigskip
\noindent Authors' addresses:
\medskip
\noindent
A. B. and R. W.: Department of Mathematics, Harvard University,
Cambridge, MA 02138 (email: bertram@math.harvard.edu,
raw@math.harvard.edu).

\medskip\noindent
G. D.:  Department of Mathematics, Princeton University,
Princeton, NJ 08544 (email: daskal@math.princeton.edu).
\end